\documentclass{article}
\usepackage[utf8]{inputenc}
\usepackage[T1]{fontenc} 
\usepackage{amsmath}
\usepackage{amsfonts}
\usepackage{amssymb}
\usepackage{amsthm}
\usepackage{color}
\usepackage{esint}
\usepackage{enumerate}
\usepackage{dsfont}
\usepackage{psfrag}
\usepackage{stmaryrd}
\usepackage{hyperref}
\usepackage{graphicx}
\usepackage{graphics}
\usepackage{epstopdf}
\usepackage{epsfig}
\usepackage{setspace}
\usepackage{tikz}
\usetikzlibrary{positioning,arrows}
\usetikzlibrary{calc,arrows}
\usetikzlibrary{er,positioning,bayesnet}

\usepackage{pgfplots}
\usepackage{pgfplotstable}
\pgfplotsset{compat=1.11} 

\usepackage{bm}
\usepackage{booktabs} 
\usepackage{caption} 
\usepackage{subcaption} 
\usepackage[all]{nowidow}
\usepackage{comment} 
\usepackage{multicol}
\usepackage{algpseudocode,algorithm,algorithmicx}
\usepackage{hyperref}
\usepackage[inline]{enumitem} 

\newtheorem{theorem}{Theorem}
\newtheorem{lemma}[theorem]{Lemma}
\newtheorem{remark}[theorem]{Remark}
\newtheorem{algo}{Algorithm}

\newcommand{\RR}{\mathbb{R}}
\newcommand{\ZZ}{\mathbb{Z}}
\newcommand{\NN}{\mathbb{N}}
\newcommand{\EE}{\mathbb{E}}
\newcommand{\eps}{\varepsilon}
\newcommand{\Id}{\text{Id}}
\newcommand{\dps}{\displaystyle}
\newcommand{\dive}{\operatorname{div}}

\definecolor{blue}{HTML}{1F77B4}
\definecolor{orange}{HTML}{FF7F0E}
\definecolor{green}{HTML}{2CA02C}

\begin{document}
\date{\today}
\author{Olga Gorynina, Claude Le Bris and Frédéric Legoll
\\
{\footnotesize \'Ecole des Ponts \& Inria}\\
{\footnotesize 6 et 8 avenue Blaise Pascal, 77455 Marne-La-Vall\'ee Cedex 2, France}\\
{\footnotesize \tt \{olga.gorynina,claude.le-bris,frederic.legoll\}@enpc.fr}\\
}
\title{Some remarks on a coupling method for the practical computation of homogenized coefficients}
\maketitle



\begin{abstract}
We numerically investigate, and improve upon, a computational approach originally introduced in~\cite{cottereau} which aims at evaluating the effective coefficient of a medium modelled by a highly oscillatory coefficient. This computational approach is based on a Arlequin type coupling. It combines the original fine-scale description of the medium (modelled by an oscillatory coefficient) with an effective description (modelled by a constant coefficient) and optimizes upon the coefficient of the effective medium to best fit the response of the actual heterogeneous medium using a purely homogeneous medium. We present here a mathematical formalization of the approach along with various improvements of the algorithms, in order to obtain a procedure as efficient as possible. Representative numerical results demonstrate the added value of our approach in comparison to the original approach.
\end{abstract}

\section{Introduction}

In this article, we investigate and improve upon a numerical approach originally introduced in~\cite{cottereau} for the practical computation of the homogenized coefficients associated to elliptic problems with highly oscillatory coefficients. The approach aims to define and construct a non-oscillating coefficient $\overline{k}$ (think e.g. of a constant matrix coefficient) that is consistent, in a sense to be made precise, with the behavior of an heterogeneous material modelled by an (possibly matrix-valued) highly oscillatory coefficient $k_\eps$. Such an approach can be considered as an alternative pathway to standard homogenization techniques when these latter are difficult to use in practice, in particular when information is missing on the coefficient $k_\eps$. Of course, modelling the material with the constant coefficient $\overline{k}$ (rather than with the function $k_\eps$) allows for much more affordable approaches, since there is no small lengthscale in that coefficient. A central question, intimately related to homogenization theory, is to construct in an efficient manner this constant matrix coefficient.

The last two authors have already addressed this question in~\cite{cocv}. The idea there was to compare the response of the actual media to various solicitations (that is, right-hand sides) with, for the same solicitation, the response of an homogeneous medium with an arbitrary constant, homogeneous coefficient. This constant coefficient was optimized upon, such that the comparison is eventually perfect.

\medskip

This article presents a mathematical and numerical study of an alternative computational approach initially introduced in~\cite{cottereau}. This approach combines the original fine-scale description of the medium (modelled by the oscillatory coefficient $k_\eps$) with an effective description (modelled by the constant coefficient $\overline{k}$). In~\cite{cottereau}, the coupling between these two models is implemented using the Arlequin technique, initially introduced in~\cite{bendhia98,arlequin} as a seamless strategy to couple two models different in nature. The specificity of the Arlequin approach (with respect to other coupling approaches) is that it smears out the interface between the models, typically as in ``$\alpha_1$ times {\it the first model} + $\alpha_2$ times {\it the second model}'', where $\alpha_1$ and $\alpha_2 = 1-\alpha_1$ are weighting functions that implicitly define the regions of the computational domain in which either of the models is used. In~\cite{cottereau}, the Arlequin approach is employed with the specific objective of practically computing homogenized coefficients. In that case, the two models coupled are the reference heterogeneous model and an arbitrary homogeneous model. The work~\cite{cottereau} proceeds by optimizing upon the coefficient of the effective medium, in order to best fit the response of the actual heterogeneous medium using a purely homogeneous medium. In the limit of asymptotically infinitely fine structures, the approach yields an approximate value of the homogenized coefficient.

We emphasize several interesting properties of the approach introduced in~\cite{cottereau}. First, this approach does not need (in sharp contrast to classical homogenization approaches) a complete knowledge of the oscillatory coefficient. It only requires to be able to solve the coupled problem defined by the Arlequin scheme. In spirit, we could even think of replacing the computation in the region modelled by $k_\eps$ by an actual experiment. Second, and again in sharp contrast to classical homogenization approaches, an approximation of $k^\star$ is obtained without computing any corrector function and without any geometric assumption on $k_\eps$ (such as periodicity, \dots) besides the fact that its homogenized limit is a constant coefficient (for the sake of simplicity, we will only consider this case hereafter; note however that the approach can be extended to cases where $k^\star$ is not constant, by suitably enlarging the search space for the tentative coefficient $\overline{k}$).

Motivated by these interesting and unusual features, our aim in this article is to mathematically formalize the problem and to investigate how to improve on the practical algorithm, in order to obtain a more efficient procedure. We refer to the upcoming companion article~\cite{ref_olga_theo} for a theoretical study of the problem.

\medskip

We note that, although different in terms of the objects that are manipulated, our former work~\cite{cocv} and the present one share many similarities in their spirit. They both are devoted to engineering-type approaches to approximate the effective coefficient of a heterogeneous medium without explicitly computing the usual ingredients of homogenization (corrector function, homogenized matrix), for example because the latter require the complete knowledge of the microstructure everywhere within the medium. In addition, they both are based on the concurrent use of the reference heterogeneous model and an arbitrary homogeneous model, and on an optimization procedure upon the coefficient of the effective medium. They however differ in how the two models are simultaneously used, and in the criterion that is minimized in the optimization loop. 

\medskip

This article is organized as follows. First, in Section~\ref{sec:state_art}, we briefly describe the Arlequin coupling method and we discuss the approach of~\cite{cottereau} in details. In Section~\ref{sec:nous}, we present several algorithmic improvements which are all illustrated with representative numerical experiments on several test cases, thereby demonstrating the added value of our approach in comparison to the original approach. In short, these improvements are the following:
\begin{itemize}
\item the Arlequin approach is based on the minimization of some energy under some constraint, and therefore leads to the introduction of a Lagrange multiplier. In Section~\ref{sec:enrich}, we suggest an enriched discretization space for this Lagrange multiplier which is better adapted to the problem at hand than the generic choice, the latter leading to a non-consistent approach. At almost no additional computational cost, this leads to a significant reduction of the error in the approximation of the homogenized coefficient. On the example we consider, the error decreases from 10\% to 0.4\%.
\item the approach of~\cite{cottereau} proceeds by optimizing upon the coefficient of the effective medium. In Section~\ref{sec:IG}, we propose a strategy to define an adequate initial guess to initialize the optimization loop, which allows the Newton algorithm to converge much faster (typically in two iterations instead of nine). 
\item in the random setting, we show how variance reduction type methods can be implemented within the approach in order to reduce the statistical noise in the approximation of the homogenized coefficient (see Section~\ref{sec:sqs}). The method suggested here, and borrowed from one of our previous work~\cite{sqs} in a different context, amounts to better selecting the realizations of the random medium. As for the other improvements, there is no additional cost in implementing this method.
\item as is well-known, the standard homogenization approach yields not only the homogenized coefficient, but also corrector functions, which are useful e.g. to reconstruct an $H^1$ approximation of the oscillatory solution. In Section~\ref{sec:corrector}, we show how to extend the approach of~\cite{cottereau} in order to reconstruct an approximation of the correctors again at no additional cost.
\end{itemize}
Several variants of the standard setting we have just described are briefly mentioned in the course of Sections~\ref{sec:state_art} and~\ref{sec:nous}. They lead to qualitatively similar theoretical conclusions and algorithmic improvements. All such variants are collected in Section~\ref{sec:variants}.

\section{Computation of the homogenized coefficient using the Arlequin approach} \label{sec:state_art}

\begin{figure}[htbp]
\centering
\includegraphics[width=0.6\textwidth]{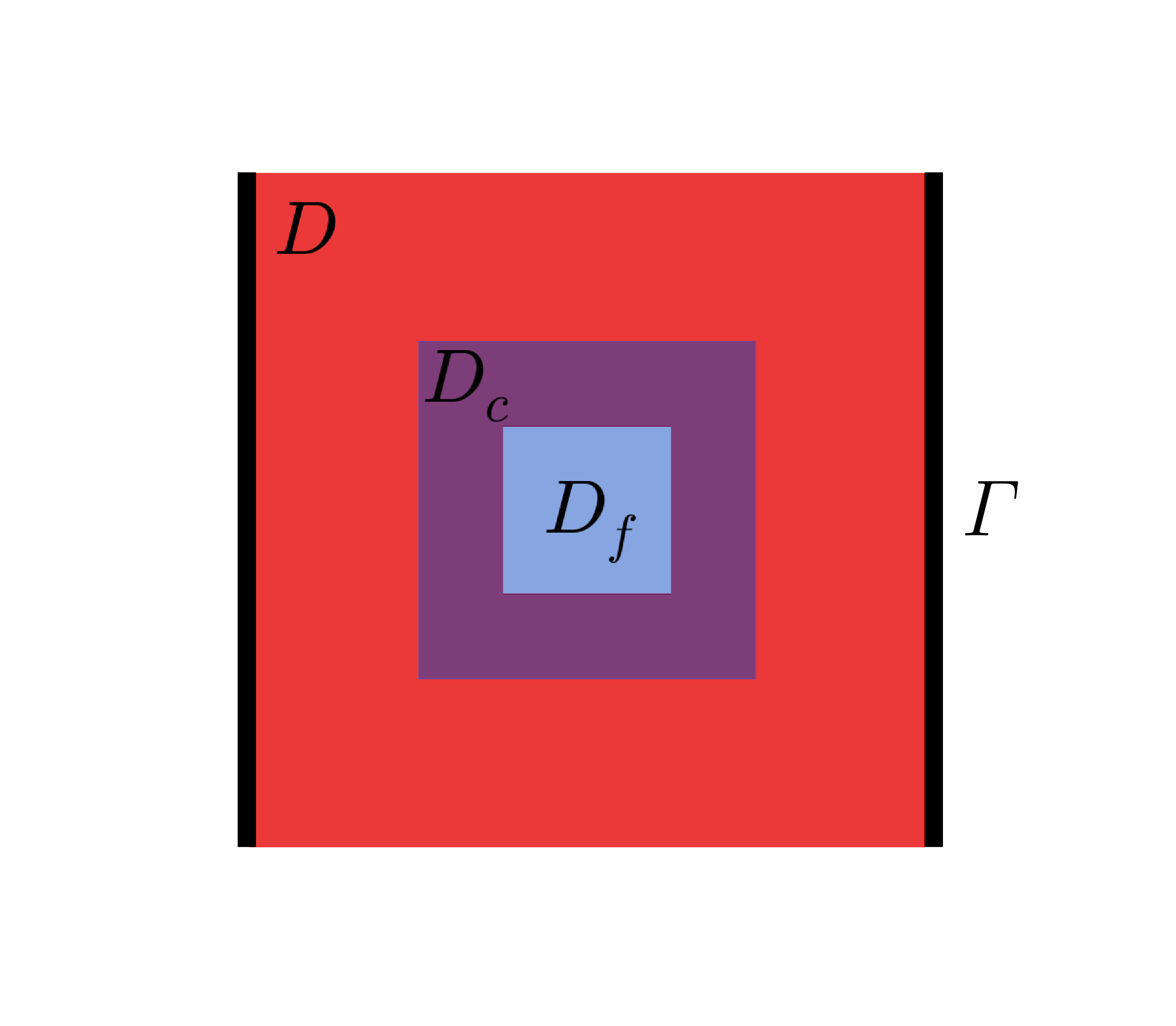}
\caption{Decomposition of the computational domain into three disjoint subdomains: a subdomain $D$ where only the effective model is defined, a subdomain $D_f$ where only the fine model is defined and a subdomain $D_c$ where both models are defined and over which they are coupled (the subscripts $f$ and $c$ obviously stand for ``fine'' and ``coupled''). Dirichlet boundary conditions are imposed on $\Gamma$, which is the part of the exterior boundary of $D$ consisting of the two black thick vertical lines. \label{fig:decompo_left}}
\bigskip
\end{figure}

The works~\cite{bendhia98,arlequin} introduce a now well-established coupling approach, namely the Arlequin approach, which (among many other applications, see~\cite{bendhia98,arlequin,cottereau_sto,rateau} for general references) is next used in~\cite{cottereau} for the specific purpose of approximating the homogenized coefficient of a heterogeneous medium. In short, the approach of~\cite{cottereau} consists in dividing the computational domain in three disjoint regions (see Figure~\ref{fig:decompo_left}). The first, inner region $D_f$ explicitly accounts for the fine-scale structure (modelled by an oscillatory, possibly matrix-valued, coefficient $k_\eps$). This first region is surrounded by a second region $D_c$, where the two models are simultaneously considered: the fine-scale structure, and the effective medium (modelled by a constant, again possibly matrix-valued, coefficient $\overline{k}$). In that second region, the two models are coupled so that, in a sense made precise in the sequel, they are consistent {\em on average}. The second region is surrounded by a third region $D$, where only the effective medium is considered, and on the exterior boundary of which suitable boundary conditions are imposed. More precisely, Dirichlet boundary conditions are imposed on some part $\Gamma$ of the exterior boundary of $D$, as shown on Figure~\ref{fig:decompo_left}. Of course, $D$, $D_c$ and $D_f$ could be arranged in a different manner, Dirichlet boundary conditions could be imposed on some other part of the exterior boundary of $D$, \dots 

For simplicity of exposition, we assume throughout the article, unless otherwise stated, that $k_\eps$, its homogenized coefficient $k^\star$ and the tentative homogeneous coefficient $\overline{k}$ are all scalar quantities. This assumption can be relaxed, as shown in Section~\ref{sec:random_nous}: in the numerical computations described there, we consider $\overline{k}$ to be a true matrix, thereby showing that the approach carries over to cases more general than the scalar case. For the well-posedness of the mathematical problem, we also assume that, on $D_c \cup D_f$, the coefficient $k_\eps$ is bounded and bounded away from 0 uniformly in $\eps$.

\medskip

The idea in~\cite{cottereau} is to optimize upon the coefficient $\overline{k}$ of the effective medium in order to best fit a linear field throughout the domain. The rationale behind this strategy is the following. Formally, the heterogeneous coefficient $k_\eps$ can be replaced by its homogenized limit $k^\star$. It is then clear that, if $\overline{k} = k^\star$, then the response of the material is linear, because the whole domain is modelled by a constant coefficient. As shown below, the converse is also true: if the response of the material is linear, then the material is homogeneous and $\overline{k} = k^\star$. Optimizing upon the coefficient $\overline{k}$ in order to best fit a linear field is thus a way to enforce $\overline{k} = k^\star$. The approach of~\cite{cottereau} thus provides a mean, alternate to standard homogenization strategies, to compute the homogenized coefficient.

\medskip

As pointed out at the end of the Introduction, the heterogeneous model and the tentative effective model can actually be coupled in various ways. We postpone a brief description of these variants until Section~\ref{sec:variants}.

\subsection{The Arlequin coupling} \label{sec:arlequin}

As mentioned above, the coupling between the two models in~\cite{cottereau} is performed using the Arlequin method, which consists in considering the following minimization problem:
\begin{equation} \label{eq:pb_min}
\inf \left\{ \begin{array}{c} {\cal E}(\overline{u},u_\eps), \quad \overline{u} \in H^1(D \cup D_c), \quad \overline{u}(x) = x_1 \ \text{on $\Gamma$}, \\ u_\eps \in H^1(D_c \cup D_f), \quad C(\overline{u}-u_\eps,\phi) = 0 \ \text{for any $\phi \in H^1(D_c)$} \end{array} \right\},
\end{equation}
where the energy ${\cal E}$ is the sum of the contributions of each of the three subdomains:
\begin{multline} \label{eq:def_E_pre}
{\cal E}(\overline{u},u_\eps) = \frac{1}{2} \int_D \overline{k} \, |\nabla \overline{u}(x)|^2 + \frac{1}{2} \int_{D_f} k_\eps(x) \, |\nabla u_\eps(x)|^2 \\ + \frac{1}{2} \int_{D_c} \Big( \frac{1}{2} \, \overline{k} \, |\nabla \overline{u}(x)|^2 + \frac{1}{2} \, k_\eps(x) \, |\nabla u_\eps(x)|^2 \Big).
\end{multline}
The last term in ${\cal E}$ accounts for the energy in the domain $D_c$, where the two models co-exist and are equally weighted. Introducing the function $\alpha_1$ defined by
\begin{equation} \label{eq:def_alpha}
\alpha_1(x) = 1 \ \ \text{in $D$},
\qquad
\alpha_1(x) = \frac{1}{2} \ \ \text{in $D_c$},
\qquad
\alpha_1(x) = 0 \ \ \text{in $D_f$},
\end{equation}
and the function $\alpha_2$ defined by $\alpha_2(x) = 1 - \alpha_1(x)$ in $D \cup D_c \cup D_f$, we can write the energy~\eqref{eq:def_E_pre} in the form
$$
{\cal E}(\overline{u},u_\eps) = \frac{1}{2} \int_{D \cup D_c \cup D_f} \alpha_1(x) \, \overline{k} \, |\nabla \overline{u}(x)|^2 + \alpha_2(x) \, k_\eps(x) \, |\nabla u_\eps(x)|^2.
$$
Note that other choices for the weighting functions $\alpha_1$ and $\alpha_2$ are possible, as discussed in Section~\ref{sec:variants_weighting} below.

\medskip

As pointed out above, we note that linear boundary conditions are enforced on $\overline{u}$ (on $\Gamma \subset \partial(D \cup D_c \cup D_f)$) in~\eqref{eq:pb_min}. Furthermore, the coarse and the fine solutions $\overline{u}$ and $u_\eps$ are consistent with one another in $D_c$, given the constraint $C(\overline{u}-u_\eps,\phi) = 0$ in~\eqref{eq:pb_min}, where $C$ is defined by
\begin{equation} \label{eq:def_C}
\forall u \in H^1(D_c), \quad \forall \phi \in H^1(D_c), \quad C(u,\phi) = \int_{D_c} \nabla u \cdot \nabla \phi + u \, \phi.
\end{equation}

It is easy to show (by considering minimizing sequences and the strong convexity of ${\cal E}$) that, for any $\overline{k} > 0$, Problem~\eqref{eq:pb_min} has a unique minimizer.

\subsection{Optimization upon the effective coefficient} \label{sec:optim}

We now explain how it is possible, by optimizing upon $\overline{k}$, to compute the homogenized coefficient $k^\star$ associated to the highly oscillatory function $k_\eps$, in the case when that homogenized limit $k^\star$ is a constant coefficient. We argue on the basis of the Arlequin coupling approach described in Section~\ref{sec:arlequin}, but we emphasize that a similar strategy can be introduced for other coupling approaches, as outlined in Section~\ref{sec:interface} below.

Consider temporarily the limit $\eps \to 0$. Then Problem~\eqref{eq:pb_min} is, at least formally (and we show by our analysis in~\cite{ref_olga_theo} that this is indeed the case after discretization), well approximated by its homogenized limit
\begin{equation} \label{eq:pb_min_star}
\inf \left\{ \begin{array}{c} {\cal E}^\star(\overline{u},u_0), \quad \overline{u} \in H^1(D \cup D_c), \quad \overline{u}(x) = x_1 \ \text{on $\Gamma$}, \\ u_0 \in H^1(D_c \cup D_f), \quad C(\overline{u}-u_0,\phi) = 0 \ \text{for any $\phi \in H^1(D_c)$} \end{array} \right\},
\end{equation}
where the energy ${\cal E}^\star$ reads as
\begin{multline} \label{eq:def_E_star}
{\cal E}^\star(\overline{u},u_0) = \frac{1}{2} \int_D \overline{k} \, |\nabla \overline{u}(x)|^2 + \frac{1}{2} \int_{D_f} k^\star \, |\nabla u_0(x)|^2 \\ + \frac{1}{2} \int_{D_c} \Big( \frac{1}{2} \, \overline{k} \, |\nabla \overline{u}(x)|^2 + \frac{1}{2} \, k^\star \, |\nabla u_0(x)|^2 \Big).
\end{multline}
In comparison to~\eqref{eq:def_E_pre}, we have replaced in the energy the oscillatory coefficient $k_\eps$ by its homogenized limit $k^\star$. We have the following result on the homogenized problem~\eqref{eq:pb_min_star}, the proof of which is postponed until Appendix~\ref{sec:proof_lemma_consistence}:

\begin{lemma} \label{lem:consistence}
  If $\overline{k} = k^\star$, then the solution to~\eqref{eq:pb_min_star} is $\overline{u}(x) = x_1$ in $D \cup D_c$ and $u_0(x) = x_1$ in $D_c \cup D_f$.

  Conversely, if $(\overline{u},u_0)$ is a solution to~\eqref{eq:pb_min_star} with $\overline{u}(x) = x_1$ in $D \cup D_c$, then $u_0(x) = x_1$ in $D_c \cup D_f$ and $\overline{k} = k^\star$.
\end{lemma}

Lemma~\ref{lem:consistence}, and in particular its second statement, motivates the idea to compare the solution $\overline{u}$ of~\eqref{eq:pb_min_star} with the reference solution $\overline{u}_{\rm ref}$ defined by $\overline{u}_{\rm ref}(x) = x_1$, which exactly corresponds to the case when $\overline{k} = k^\star$. Stated otherwise, we optimize upon $\overline{k}$ by considering the minimization problem
\begin{equation} \label{eq:optim_J}
\inf \left\{ J(\overline{k}), \quad \overline{k} \in (0,\infty) \right\},
\end{equation}
with
\begin{equation} \label{eq:def_J}
J(\overline{k}) = \int_{D \cup D_c} | \nabla \overline{u} - \nabla \overline{u}_{\rm ref} |^2 = \int_{D \cup D_c} | \nabla \overline{u} - e_1 |^2,
\end{equation}
where $(\overline{u},u_0)$ is the solution to~\eqref{eq:pb_min_star} (that depends on $\overline{k}$). In the homogenized limit considered here, the infimum of $J$ is zero, and is attained if and only if $\overline{k} = k^\star$.

\medskip

In practice, the parameter $\eps$ is small but we cannot consider the limit $\eps \to 0$ (and therefore the energy~\eqref{eq:def_E_star}), because we do not know $k^\star$ beforehand (the whole point of the approach we are describing being of course to \emph{determine} this homogenized coefficient). In addition, Problems~\eqref{eq:pb_min_star} and~\eqref{eq:pb_min} can of course not be exactly solved in practice, but should be discretized.

A natural method is then to proceed as follows: we minimize~\eqref{eq:def_J} with respect to the coefficient $\overline{k}$, where $(\overline{u},u_\eps)$ is the solution to~\eqref{eq:pb_min} (or the solution to its discretized version that we present in Section~\ref{sec:disc} below, see problem~\eqref{eq:pb_min_H}). We have the following results:
\begin{itemize}
\item as shown in~\cite{ref_olga_theo}, when considering the discretized version of~\eqref{eq:pb_min} ($\overline{u}$ is then discretized on a mesh of size $H$ independent of $\eps$ and $u_\eps$ is discretized on a mesh of size $h \ll \eps$), the minimization procedure~\eqref{eq:optim_J}--\eqref{eq:def_J} yields an optimal value $\overline{k}^{\rm opt}(\eps,h,H)$ which itself converges to $k^\star$ when $\eps$, $h$ and $H$ tend to 0. More precisely, we show $\dps \lim_{H \to 0} \lim_{\eps \to 0} \lim_{h \to 0} \overline{k}^{\rm opt}(\eps,h,H) = k^\star$.
\item in the absence of any discretization, the approach {\em does not} yield the value of homogenized coefficient: if we minimize~\eqref{eq:def_J} with respect to $\overline{k}$, where $(\overline{u},u_\eps)$ is the solution to the {\em exact} problem~\eqref{eq:pb_min}, then we find an optimal value $\overline{k}^{\rm opt}(\eps)$ which is {\em different} from the homogenized coefficient $k^\star$ we seek, even after passing to the limit $\eps \to 0$. This issue is related to the fact that, in the absence of any discretization, the constraint $C(\overline{u}-u_\eps,\phi) = 0$ for any $\phi \in H^1(D_c)$ in~\eqref{eq:pb_min} imposes $\overline{u} = u_\eps$ in $D_c$, which is a strong requirement. The problem~\eqref{eq:pb_min} is thus equivalent to minimizing
\begin{equation} \label{eq:nrj_sigma}
\frac{1}{2} \int_D \overline{k} \, |\nabla \overline{u}(x)|^2 + \frac{1}{2} \int_{D_c \cup D_f} \sigma_\eps(x) \, |\nabla u_\eps(x)|^2,
\end{equation}
with the transmission condition $\overline{u} = u_\eps$ across the interface between $D$ and $D_c \cup D_f$ (and of course the boundary condition for $\overline{u}$ on $\Gamma$), where the oscillatory coefficient $\sigma_\eps$ is defined by $\dps \sigma_\eps(x) = \frac{\overline{k}+k_\eps(x)}{2}$ in $D_c$ and $\sigma_\eps(x) = k_\eps(x)$ in $D_f$. The homogenized limit of~\eqref{eq:nrj_sigma} is then
$$
\frac{1}{2} \int_D \overline{k} \, |\nabla \overline{u}(x)|^2 + \frac{1}{2} \int_{D_c \cup D_f} \sigma^\star \, |\nabla u_0(x)|^2,
$$
with $\dps \sigma^\star = \frac{(\overline{k}+k)^\star}{2}$ in $D_c$ and $\sigma^\star = k^\star$ in $D_f$. The homogenized coefficient $(\overline{k}+k)^\star$ is in general different from $\overline{k}+k^\star$, thus the fact that the optimal value $\overline{k}^{\rm opt}(\eps)$ converges to a value different from $k^\star$. The above inconsistency is even more evident in dimension $d=1$, where it is possible to analytically solve~\eqref{eq:pb_min}: assuming that $k_\eps = k_{\rm per}(\cdot/\eps)$ for a $\ZZ$-periodic function $k_{\rm per}$ and that the size of $D_c$ and $D_f$ is a multiple of $\eps$, we find that $\overline{u}'(x) = m/\overline{k}$ in $D$, where $m$ is independent of $\eps$, depends on $\overline{k}$ and satisfies
$$
| D \cup D_c \cup D_f | = m \left( \frac{|D|}{\overline{k}} + \frac{|D_f|}{k^\star} + |D_c| \, \left\langle \frac{2}{\overline{k} + k_{\rm per}} \right\rangle \right),
$$
where $\langle \cdot \rangle$ denotes the average over the periodic cell. Considering (for the sake of simplicity) the minimization of $\dps J_{\rm 1D}(\overline{k}) = \int_D | \overline{u}' - 1|^2$ instead of~\eqref{eq:def_J} yields the optimal value $\overline{k}^{\rm opt}$ which is defined as the unique solution to
$$
| D_c \cup D_f | = |D_f| \ \frac{\overline{k}^{\rm opt}}{k^\star} + |D_c| \ \overline{k}^{\rm opt} \, \left\langle \frac{2}{\overline{k}^{\rm opt} + k_{\rm per}} \right\rangle.
$$
Using Jensen inequality, we observe that $\overline{k}^{\rm opt} > k^\star$ as soon as $k_{\rm per}$ is not constant. A similar conclusion is reached when using the actual function $J$ rather than $J_{\rm 1D}$: $\overline{k}^{\rm opt} \neq k^\star$. 
\end{itemize}
Considering the discretized version of~\eqref{eq:pb_min} is thus critical, not only for practical reasons, but also to have a method that converges to the correct limit. We note that this is not the case for the homogenized problem~\eqref{eq:pb_min_star}, as pointed out above using Lemma~\ref{lem:consistence}.

\subsection{Discretization} \label{sec:disc}

\begin{figure}[htbp]
\centering
\includegraphics[width=0.6\textwidth]{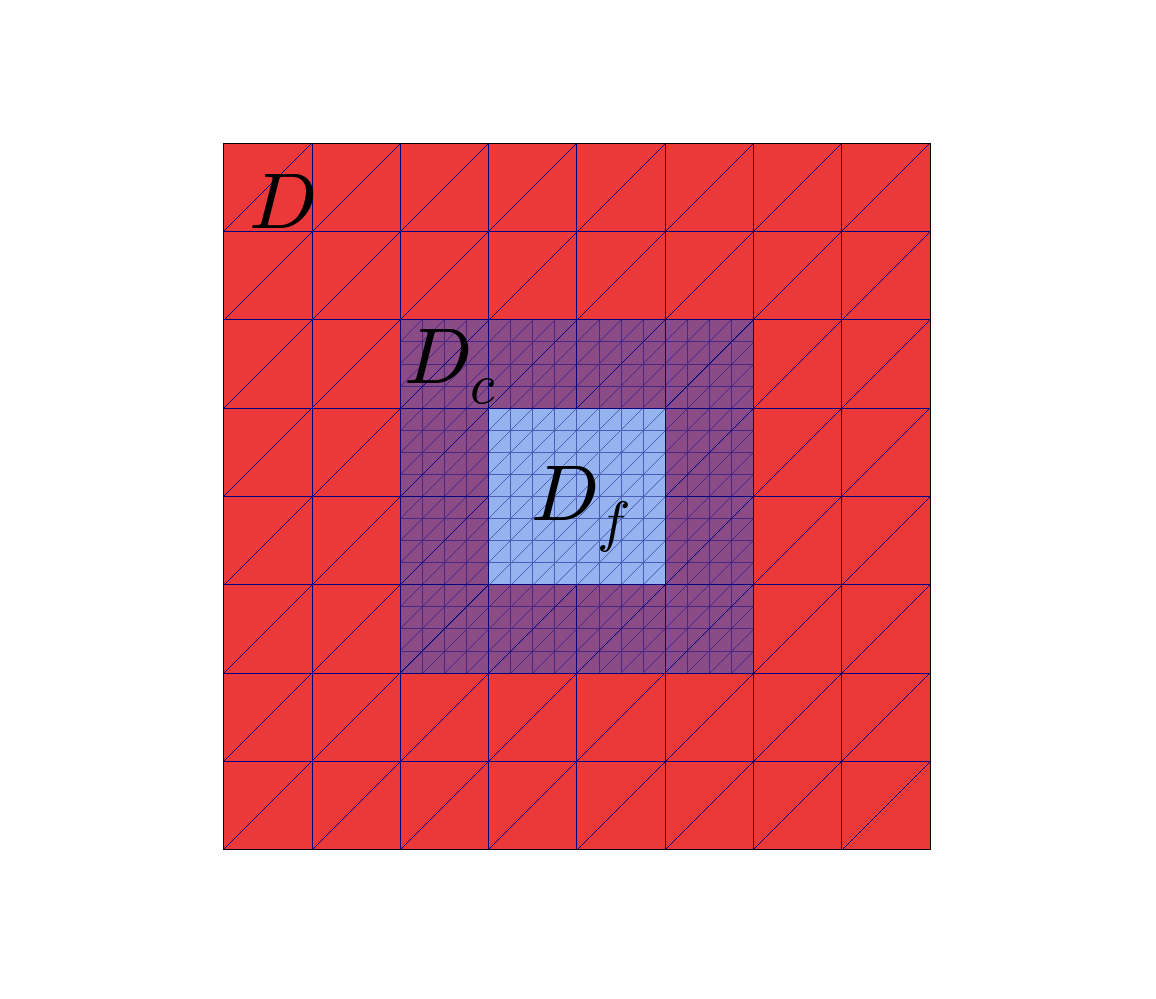}
\caption{A coarse (resp. fine) mesh is used in $D \cup D_c$ (resp. $D_c \cup D_f$). \label{fig:decompo_right}}
\bigskip
\end{figure}

We introduce a coarse mesh ${\cal T}_H$ (of mesh size $H$) in the subdomains $D$ and $D_c$ and a fine mesh ${\cal T}_h$ (of mesh size $h$) in the subdomains $D_c$ and $D_f$ (see Figure~\ref{fig:decompo_right}). We assume that the coarse meshes of $D$ and $D_c$ are consistent with one another on $\partial D \cap \partial D_c$, namely that they match on the interface (and likewise for the fine meshes of $D_c$ and $D_f$).
We also assume that, in $D_c$, the fine mesh is a submesh of the coarse mesh. The fine mesh size $h$ is assumed to be much smaller than the finest oscillations of the coefficient $k_\eps$: $h \ll \eps$ (say $h \approx \eps/10$). In contrast, the coarse mesh size $H$ can be chosen independent of $\eps$, and therefore satisfies $H \gg h$. We next introduce the corresponding finite element spaces:
\begin{gather}
  V_H =  \left\{ u \in H^1(D \cup D_c), \quad \text{$u$ is piecewise affine on the coarse mesh ${\cal T}_H$} \right\},
  \nonumber
  \\
  V_h = \left\{ u \in H^1(D_c \cup D_f), \quad \text{$u$ is piecewise affine on the fine mesh ${\cal T}_h$} \right\},
  \nonumber
  \\
  W_H =  \left\{ \phi \in H^1(D_c), \quad \text{$\phi$ is piecewise affine on the coarse mesh ${\cal T}_H$} \right\}.
  \label{eq:def_WH}
\end{gather}

\medskip

\noindent
The minimization problem~\eqref{eq:pb_min} is then approximated by
\begin{equation} \label{eq:pb_min_H}
\inf \left\{ \begin{array}{c} {\cal E}(\overline{u},u_\eps), \quad \overline{u} \in V_H, \quad \overline{u}(x) = x_1 \ \text{on $\Gamma$}, \\ u_\eps \in V_h, \quad C(\overline{u}-u_\eps,\phi) = 0 \ \text{for any $\phi \in W_H$} \end{array} \right\},
\end{equation}
where the energy ${\cal E}$ and the constraint $C$ are defined by~\eqref{eq:def_E_pre} and~\eqref{eq:def_C}. Similarly to~\eqref{eq:pb_min}, problem~\eqref{eq:pb_min_H} has a unique minimizer (considering again minimizing sequences and using the strong convexity of ${\cal E}$).

We observe that the constraint $C(\overline{u}-u_\eps,\phi) = 0$ for any $\phi \in W_H$ means that, on $D_c$, $\overline{u}$ is the projection on piecewise affine functions on the coarse mesh ${\cal T}_H$ of $u_\eps$, itself a piecewise affine function on ${\cal T}_h$. In that sense, on $D_c$, $\overline{u}$ and $u_\eps$ are consistent with one another {\em on average}. This is a much looser constraint than in the continuous context, where the constraint implied $\overline{u} = u_\eps$ on $D_c$ (see~\eqref{eq:pb_min} and~\eqref{eq:def_C}).

By convexity of the energy and linearity of the constraint (see below), solving the minimization problem~\eqref{eq:pb_min_H} is equivalent to solving the following variational formulation: find $\overline{u} \in V_H^{\rm Dir BC}$, $u_\eps \in V_h$ and $\psi \in W_H$ such that
\begin{equation} \label{eq:arlequin}
  \begin{cases}
  \forall \overline{v} \in V_H^0, \quad &\overline{A}(\overline{u},\overline{v}) + C(\overline{v},\psi) = 0,
  \\
  \forall v_\eps \in V_h, \quad & A_\eps(u_\eps,v_\eps) - C(v_\eps,\psi) = 0,
  \\
  \forall \phi \in W_H, \quad & C(\overline{u}-u_\eps,\phi) = 0,
  \end{cases}
\end{equation}
where
\begin{gather*}
V_H^{\rm Dir BC} = \left\{ v \in V_H, \quad \text{$v(x)=x_1$ on $\Gamma$} \right\},
\\
V_H^0 = \left\{ v \in V_H, \quad \text{$v=0$ on $\Gamma$} \right\},
\end{gather*}
and where the bilinear forms $\overline{A}$ and $A_\eps$ are respectively defined by
\begin{equation}
  \overline{A}(\overline{u},\overline{v})
  =
  \int_D \overline{k} \, \nabla \overline{u}(x) \cdot \nabla \overline{v}(x) + \frac{1}{2} \int_{D_c} \overline{k} \, \nabla \overline{u}(x) \cdot \nabla \overline{v}(x) 
  \label{eq:def_Abar}
\end{equation}
and
\begin{equation}
  A_\eps(u_\eps,v_\eps)
  =
  \frac{1}{2} \int_{D_c} k_\eps(x) \, \nabla u_\eps(x) \cdot \nabla v_\eps(x) + \int_{D_f} k_\eps(x) \, \nabla u_\eps(x) \cdot \nabla v_\eps(x).
  \label{eq:def_Aeps}
\end{equation}
Indeed, any minimizer of~\eqref{eq:pb_min_H} is a solution to the associated Euler-Lagrange equations, which read as~\eqref{eq:arlequin}. Conversely, a direct computation shows that, if $(\overline{u},u_\eps,\psi) \in V_H^{\rm Dir BC} \times V_h \times W_H$ is a solution to~\eqref{eq:arlequin}, then $(\overline{u},u_\eps)$ is a minimizer of~\eqref{eq:pb_min_H}.

We conclude this section by observing that~\eqref{eq:arlequin} has a unique solution for any $\overline{k} > 0$. Consider indeed two solutions $(\overline{u}^j,u_\eps^j,\psi^j)$, $j=1,2$, of~\eqref{eq:arlequin}. We then obtain that $(\overline{u}^j,u_\eps^j)$ is a minimizer of~\eqref{eq:pb_min_H}, and therefore $\overline{u}^1 = \overline{u}^2$ and $u_\eps^1 = u_\eps^2$. The second line of~\eqref{eq:arlequin} hence yields
$$
\forall v_\eps \in V_h, \quad C(v_\eps,\psi^1-\psi^2) = 0.
$$
We can take $v_\eps \in V_h$ such that, on $D_c$, $v_\eps = \psi^1-\psi^2$ (recall that the fine mesh in $D_c$ is a submesh of the coarse mesh). We therefore obtain $C(\psi^1-\psi^2,\psi^1-\psi^2) = 0$, which yields $\psi^1 = \psi^2$ and therefore the uniqueness of the solution to~\eqref{eq:arlequin}.

\subsection{Extension of the approach to random cases} \label{sec:random_cottereau}

We now consider the case of materials modelled by a {\em random} oscillatory coefficient. Identifying the homogenized coefficient of such materials is known to be a challenging problem from the practical viewpoint. The approach introduced in~\cite{cottereau} was actually used there for this specific purpose. With the aim to couple the material modelled by $\overline{k}$ to an {\em average} (over the random realizations) of the heterogeneous materials modelled by $k_\eps(\cdot,\omega)$, we now revisit the approach described above.

We assume that the fine-scale model is described by the highly oscillatory random scalar function $k_\eps(x,\omega)$. We also assume that the homogenized limit $k^\star$ associated to $k_\eps$ is a constant and deterministic coefficient (this is for instance the case when $k_\eps(x,\omega) = k_{\rm sto}(x/\eps,\omega)$ for a stationary and ergodic coefficient $k_{\rm sto}$), and that $k^\star$ is scalar-valued. In the random setting (see~\cite{cottereau_sto}), the Arlequin approach consists in considering the minimization problem
\begin{equation} \label{eq:pb_min_sto}
\inf \left\{ \begin{array}{c} {\cal E}_{\rm sto}(\overline{u},u_\eps), \quad \overline{u} \in H^1(D \cup D_c), \quad \overline{u}(x) = x_1 \ \text{on $\Gamma$}, \\ u_\eps \in V_f, \quad C_{\rm sto}(\overline{u}-u_\eps,\phi) = 0 \ \text{for any $\phi \in W_{\rm sto}$} \end{array} \right\},
\end{equation}
where $\dps V_f = L^2\big(\Omega,H^1(D_c \cup D_f)\big)$, $\Omega$ is the probability space and where the energy ${\cal E}_{\rm sto}$ is
\begin{multline*} 
{\cal E}_{\rm sto}(\overline{u},u_\eps) = \frac{1}{2} \int_D \overline{k} \, |\nabla \overline{u}(x)|^2 + \frac{1}{2} \int_{D_f} \EE \Big[ k_\eps(x,\cdot) \, |\nabla u_\eps(x,\cdot)|^2 \Big] \\ + \frac{1}{2} \int_{D_c} \Big( \frac{1}{2} \, \overline{k} \, |\nabla \overline{u}(x)|^2 + \frac{1}{2} \, \EE \Big[ k_\eps(x,\cdot) \, |\nabla u_\eps(x,\cdot)|^2 \Big] \Big).
\end{multline*}
We note that the function $\overline{u}$ is deterministic while the function $u_\eps$ is random. Again, the coefficient $\overline{k}$ is constant and deterministic, which is reminiscent of the fact that $k^\star$ is constant and deterministic. In the energy ${\cal E}_{\rm sto}$, we consider the expectation of the random contribution. The constraint also involves an average over the realizations: it reads as
$$
C_{\rm sto}(u,\phi) = \EE \left[ \int_{D_c} \nabla u \cdot \nabla \phi + u \, \phi \right].
$$
The Lagrange multiplier space (which used to simply be $H^1(D_c)$ in the deterministic context, see~\eqref{eq:pb_min}) is now given by
\begin{equation} \label{eq:W_sto}
  W_{\rm sto}
  =
  \big\{ \psi + \theta, \quad \psi \in H^1_{\rm mean}(D_c), \quad \theta \in L^2(\Omega) \big\},
\end{equation}
where
$$
H^1_{\rm mean}(D_c) = \left\{ \psi \in H^1(D_c), \quad \int_{D_c} \psi = 0 \right\}.
$$
Similarly to Problem~\eqref{eq:pb_min}, it is straightforward to show that, for any $\overline{k} > 0$, Problem~\eqref{eq:pb_min_sto} has a unique minimizer (we refer to~\cite{ref_olga_theo} for a proof of this claim, along with some justifications for the above choices).

\medskip

Solving the minimization problem~\eqref{eq:pb_min_sto} is equivalent to solving its associated Euler-Lagrange equations, which read as follows: find $\overline{u} \in H^1(D \cup D_c)$ with $\overline{u}(x) = x_1$ on $\Gamma$, $u_\eps \in V_f$ and $\psi_{\rm sto} \in W_{\rm sto}$ such that 
\begin{equation} \label{eq:arlequin_continuous_sto}
  \begin{cases}
  \forall \overline{v} \in V^0, & \quad \overline{A}(\overline{u},\overline{v}) + C_{\rm sto}\left(\overline{v},\psi_{\rm sto}\right) = 0,
  \\
  \forall v_\eps \in V_f, & \quad A_\eps(u_\eps,v_\eps) - C_{\rm sto}(v_\eps,\psi_{\rm sto}) = 0,
  \\
  \forall \phi_{\rm sto} \in W_{\rm sto}, & \quad C_{\rm sto}(\overline{u} - u_\eps,\phi_{\rm sto}) = 0,
  \end{cases}
\end{equation}
where
$$ 
V^0 = \left\{ v \in H^1(D \cup D_c) \quad \text{such that} \quad v = 0 \ \text{on} \ \Gamma \right\},
$$
where the bilinear form $\overline{A}$ is defined by~\eqref{eq:def_Abar}, and where the bilinear form $A_\eps$ is defined, for any $u$ and $v$ in $V_f$, by (compare with~\eqref{eq:def_Aeps})
$$
A_\eps(u,v)
=
\frac{1}{2} \int_{D_c} \EE \big[ k_\eps(x,\cdot) \, \nabla u \cdot \nabla v \big] + \int_{D_f} \EE \big[ k_\eps(x,\cdot) \, \nabla u \cdot \nabla v \big].
$$
In view of~\eqref{eq:W_sto}, the functions $\psi_{\rm sto}$ and $\phi_{\rm sto}$ have the following form:
$$
\psi_{\rm sto}(x,\omega) = \psi(x) + \theta(\omega) \qquad \text{and} \qquad \phi_{\rm sto}(x,\omega) = \phi(x)+\xi(\omega),
$$
where $\psi,\phi \in H^1_{\rm mean}(D_c)$ and $\theta,\xi \in L^2(\Omega)$. The system~\eqref{eq:arlequin_continuous_sto} thus reads as: find $\overline{u} \in H^1(D \cup D_c)$ with $\overline{u}(x) = x_1$ on $\Gamma$, $u_\eps \in V_f$, $\psi \in H^1_{\rm mean}(D_c)$ and $\theta \in L^2(\Omega)$ such that 
\begin{equation} \label{eq:arlequin_original}
  \begin{cases}
  \forall \overline{v} \in V^0, & \dps \quad \overline{A}(\overline{u},\overline{v}) + C(\overline{v},\psi)+\EE \left[ \theta \int_{D_c} \overline{v}\right] = 0,
  \\
  \forall v_\eps \in V_f, & \dps \quad A_\eps(u_\eps,v_\eps) - C(\EE \left[ v_\eps \right],\psi) - \EE \left[ \theta \int_{D_c} v_\eps \right] = 0,
  \\
  \forall \phi \in H^1_{\rm mean}(D_c), & \quad C(\overline{u} - \EE \left[u_\eps \right],\phi) = 0,
  \\
  \forall \xi \in L^2(\Omega), & \dps \quad \EE \left[\xi \int_{D_c}\left(\overline{u} -u_\eps\right) \right] = 0,
  \end{cases}
\end{equation}
where we recall that the bilinear form $C$ is given by~\eqref{eq:def_C} (note that the system~\eqref{eq:arlequin_continuous_sto} involves $C_{\rm sto}$, while~\eqref{eq:arlequin_original} involves $C$). Using arguments similar to those used to show that~\eqref{eq:arlequin} is well-posed (see the end of Section~\ref{sec:disc}), we know that Problem~\eqref{eq:arlequin_original} has a unique solution.

Taking $\xi \equiv 1$ in the fourth line of~\eqref{eq:arlequin_original}, we see that the mean over $D_c$ of $\overline{u} - \EE \left[u_\eps\right]$ vanishes. The third line of~\eqref{eq:arlequin_original} then implies $\overline{u}(x) = \EE \left[u_\eps(x,\cdot)\right]$ for any $x \in D_c$. 

\medskip

Performing a Galerkin discretization of~\eqref{eq:arlequin_original} is not straightforward, because one has to introduce a finite dimensional subspace of the space $H^1_{\rm mean}(D_c)$ of functions with vanishing mean.
In practice, we prefer to recast~\eqref{eq:arlequin_original} in a slightly different formulation: find $\overline{u} \in H^1(D \cup D_c)$ with $\overline{u}(x) = x_1$ on $\Gamma$, $u_\eps \in V_f$ and $\psi \in H^1(D_c)$ such that
\begin{equation} \label{eq:arlequin_alternative}
  \begin{cases}
  \forall \overline{v} \in V^0, & \quad \overline{A}(\overline{u},\overline{v}) + C(\overline{v},\psi) = 0,
  \\
  \forall v_\eps \in H^1(D_c \cup D_f), & \quad A_\eps^\omega(u_\eps(\omega),v_\eps) - C(v_\eps,\psi) = 0 \ \ \text{a.s.},
  \\
  \forall \phi \in H^1(D_c), & \quad C(\overline{u} - \EE \left[u_\eps \right],\phi) = 0,
  \end{cases}
\end{equation}
where the bilinear form $A_\eps^\omega$ is given by
\begin{equation}
  A_\eps^\omega(u,v)
  =
  \frac{1}{2} \int_{D_c} k_\eps(x,\omega) \, \nabla u \cdot \nabla v + \int_{D_f} k_\eps(x,\omega) \, \nabla u \cdot \nabla v.
  \label{eq:def_aeps_omega}
\end{equation}
We now discuss the relation between Problem~\eqref{eq:arlequin_original} and Problem~\eqref{eq:arlequin_alternative}. We first recall that Problem~\eqref{eq:arlequin_original} is well-posed. We have the following result, the proof of which is postponed until Appendix~\ref{app:equivalence}.

\begin{lemma} \label{lem:dimanche}
  Let $(\overline{u},u_\eps,\psi,\theta)$ be the solution to~\eqref{eq:arlequin_original}. Then $(\overline{u},u_\eps,\psi)$ is a solution to~\eqref{eq:arlequin_alternative}. Conversely, let $(\widetilde{\overline{u}},\widetilde{u_\eps},\widetilde{\psi})$ be a solution to~\eqref{eq:arlequin_alternative}. Then $\widetilde{\overline{u}} = \overline{u}$, $\widetilde{\psi} = \psi$ and $\widetilde{u_\eps}(x,\omega) = u_\eps(x,\omega) + c_1(\omega)$ with $\EE[c_1] = 0$. The solution to~\eqref{eq:arlequin_alternative} is thus unique up to the addition, for $u_\eps$, of a random constant of vanishing expectation.
\end{lemma}

The solution to~\eqref{eq:arlequin_alternative} is not unique, but this does not raise any issue in our context. Indeed, our aim, following~\eqref{eq:optim_J}, is to minimize the function $J$ defined by~\eqref{eq:def_J}. Since $J(\overline{k})$ only depends on $\overline{u}$ (and in particular does not depend on $u_\eps$), we can equivalently compute $J(\overline{k})$ by solving~\eqref{eq:arlequin_original} or~\eqref{eq:arlequin_alternative}. In the vein of the discussion at the end of Section~\ref{sec:optim}, we note that we need to consider the discretized version of~\eqref{eq:arlequin_original} or~\eqref{eq:arlequin_alternative} to ensure that the optimization strategy indeed converges to $k^\star$.

\medskip

Problem~\eqref{eq:arlequin_alternative} needs to be discretized in space (through e.g. the introduction of finite element spaces) and in probability (using the Monte Carlo method: introduction of several realizations and replacement of expectations by averages over realizations).

Considering $M$ realizations, the equations~\eqref{eq:arlequin_alternative} then become: find $\overline{u} \in V_H^{\rm Dir BC}$, $u_\eps^m \in V_h$ for any $1 \leq m \leq M$ and $\psi \in W_H$ such that
\begin{equation} \label{eq:arlequin_sto}
  \begin{cases}
  \forall \overline{v} \in V_H^0, \quad & \overline{A}(\overline{u},\overline{v}) + C(\overline{v},\psi) = 0,
  \\
  \forall 1 \leq m \leq M, \ \ \forall v_\eps^m \in V_h, \quad & A_\eps^m(u_\eps^m,v_\eps^m) - C(v_\eps^m,\psi) = 0,
  \\
  \forall \phi \in W_H, \quad & \dps C(\overline{u},\phi) - \frac{1}{M} \sum_{m=1}^M C(u_\eps^m,\phi) = 0,
  \end{cases}
\end{equation}
where the finite dimensional spaces $V_H^{\rm Dir BC}$, $V_H^0$, $V_h$ and $W_H$ have been introduced in Section~\ref{sec:disc} and where the bilinear form $A_\eps^m$ is given (compare to~\eqref{eq:def_aeps_omega}) by
$$
A_\eps^m(u,v) = \frac{1}{2} \int_{D_c} k_\eps(x,\omega_m) \, \nabla u \cdot \nabla v + \int_{D_f} k_\eps(x,\omega_m) \, \nabla u \cdot \nabla v,
$$
where $k_\eps(x,\omega_m)$ is the $m$-th realization of the random function $k_\eps(x,\cdot)$.

\medskip

Problem~\eqref{eq:arlequin_sto} is huge since it involves $M$ fine-scale functions $u_\eps^m$, which are all discretized on a fine mesh of size $h \ll \eps$. Furthermore, $\overline{u}$ and all the functions $u_\eps^m$ are coupled with one another through the Lagrange multiplier $\psi$ (see the first and second lines of~\eqref{eq:arlequin_sto}). As such, problem~\eqref{eq:arlequin_sto} thus seems intractable. In Section~\ref{sec:random_nous} below, we will return to this problem and explain how to recast it in a convenient way, amenable to be solved at a limited computational cost.

\section{Our contributions} \label{sec:nous}

We now present our contributions on this problem, following the executive summary provided at the end of the Introduction.
In Section~\ref{sec:enrich}, we first suggest a choice of the approximation space for the Lagrange multiplier $\psi$ that is better adapted to the problem at hand than the generic choice~\eqref{eq:def_WH}. Second, we present a strategy for suitably initializing the optimization loop on the effective coefficient $\overline{k}$ (see Section~\ref{sec:IG}). We next describe in Sections~\ref{sec:random_nous} and~\ref{sec:sqs} how the random problem of Section~\ref{sec:random_cottereau} can be efficiently solved. We eventually consider in Section~\ref{sec:corrector} the reconstruction of the corrector function, a key ingredient in homogenization theory. 

\subsection{Enrichment of the Lagrange multiplier space} \label{sec:enrich}

As seen above, the Arlequin approach (used here for the specific aim of determining a homogenized coefficient) amounts to finding some $(\overline{u},u_\eps,\psi)$ solution to~\eqref{eq:arlequin}. For any $\eps > 0$, none of these functions (which all depend on $\eps$) has a simple expression. However, in the homogenized limit $\eps \to 0$, it turns out that the Lagrange multiplier $\psi$ has a simple expression, as we now explain.

Consider formally the limit $\eps \to 0$ again, and assume that the optimization procedure described in Section~\ref{sec:optim} has converged: we can thus suppose that $\overline{k} = k^\star$ and $\overline{u}(x) = x_1$ in $D \cup D_c$. At the continuous level, the first line of~\eqref{eq:arlequin} then reads as
$$
\forall \overline{v} \in V^0, \quad \int_D k^\star \, e_1 \cdot \nabla \overline{v} + \frac{1}{2} \int_{D_c} k^\star \, e_1 \cdot \nabla \overline{v} + C(\overline{v},\psi) = 0.
$$
The function $\psi$ is thus the solution to
$$
\begin{cases}
- \Delta \psi + \psi = 0 \quad \text{in $D_c$},
\\
\dps \nabla \psi \cdot n = \frac{k^\star}{2} \, e_1 \cdot n \quad \text{on $\partial D \cap \partial D_c$}, \quad \nabla \psi \cdot n = -\frac{k^\star}{2} \, e_1 \cdot n \ \ \text{on $\partial D_c \cap \partial D_f$}.
\end{cases}
$$
Since $k^\star$ is a scalar constant, we have $\psi = k^\star \, \psi_0$ with $\psi_0$ solution to 
\begin{equation} \label{eq:def_psi0}
\hspace{-3mm} \left\{ \begin{array}{l}
- \Delta \psi_0 + \psi_0 = 0 \quad \text{in $D_c$},
\\
\dps \nabla \psi_0 \cdot n = \frac{1}{2} \, e_1 \cdot n \ \text{on $\partial D \cap \partial D_c$}, \ \ \nabla \psi_0 \cdot n = -\frac{1}{2} \, e_1 \cdot n \ \text{on $\partial D_c \cap \partial D_f$}.
\end{array} \right.
\end{equation}
We thus observe that $\psi_0$ can be computed independently of the knowledge of $k^\star$. This motivates the following enrichment procedure. Instead of considering, in~\eqref{eq:arlequin}, the space $W_H$ of the functions of $H^1(D_c)$ that are piecewise affine on the coarse mesh ${\cal T}_H$, we consider the space
$$
W_H^{\rm enrich} = W_H + \text{Span $\psi_0$}.
$$
With this choice, the approach~\eqref{eq:arlequin} is now consistent.
We then hope (and this will indeed be the case in the numerical examples discussed below) to reduce the error with respect to $H$ of the approach. We are indeed enlarging the discretization space, so that the exact solution (at convergence $\eps \to 0$) of the problem belongs to that space. 

Of course, the solution $\psi_0$ to~\eqref{eq:def_psi0} is only approximated numerically in practice, using e.g. a finite element approach. This computation has only to be performed once, independently of the number of iterations to solve the minimization problem~\eqref{eq:optim_J}. The additional cost can thus be neglected.

Note also that $\psi_0$ should be (and can be) approximated with an accuracy much better than that obtained using a finite element method with piecewise affine functions on the coarse mesh ${\cal T}_H$ (otherwise, no accuracy gain can be expected). The mesh used to approximate $\psi_0$ should hence be of size (much) smaller than $H$.

\medskip

\begin{figure}[htbp]
\centering
\includegraphics[width=0.8\textwidth]{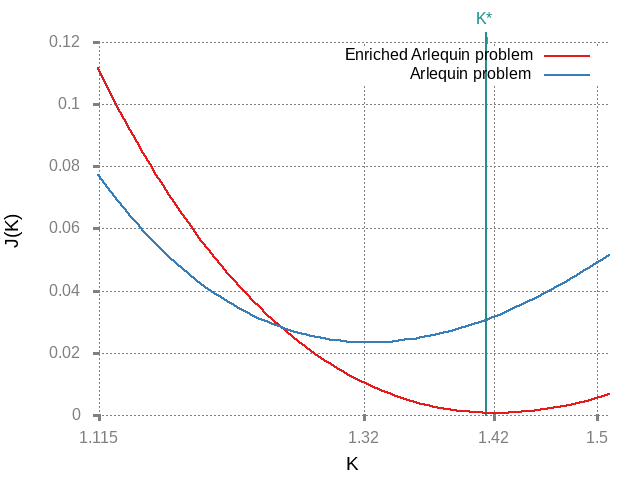}
\caption{Plot of the function $\overline{k} \mapsto J(\overline{k})$ in the case when we use $W_H^{\rm enrich}$ (red curve, which attains its minimum at 1.42) for the Lagrange multiplier space, or simply the space $W_H$ (blue curve, which attains its minimum at 1.32). The exact value $k^\star=\sqrt{2} \approx 1.4142$ is shown by the green vertical line. \label{fig:enrich}}
\bigskip
\end{figure}

As shown on Figure~\ref{fig:enrich}, the enrichment of the Lagrange multiplier space turns out to be very beneficial. The computations are performed in the following two-dimensional case. Henceforth, we consider the domains $D_f = (-1,1) \times (-1,1)$, $D_c = (-2,2) \times (-2,2) \setminus \overline{D_f}$ and $D = (-4,4) \times (-4,4) \setminus \overline{(D_c \cup D_f)}$ (see Figure~\ref{fig:decompo_left}). On $D_c \cup D_f$, we consider a periodic checkerboard geometry for the heterogeneous model: $k_\eps$ is piecewise constant on each cell of size $\eps \times \eps$, and is equal to either 1 or 2. We set $\eps = 1/16$, and use a fine mesh of size $h = 1/64$. 
The coarse mesh size is $H = 1$. In that case, the homogenized coefficient is known (although of course not used in the approach) and is equal to $k^\star = \sqrt{2}$.

We first plot (see the blue curve) the function $\overline{k} \mapsto J(\overline{k})$ defined by~\eqref{eq:def_J}, where $(\overline{u},u_\eps)$ is the solution to~\eqref{eq:pb_min_H} (or equivalently~\eqref{eq:arlequin}). We see that this function attains its minimum for $\overline{k} \approx 1.32$. Since $k^\star = \sqrt{2} \approx 1.4142$ (a value represented by the green vertical line on Figure~\ref{fig:enrich}), we observe that the error on the evaluation of the homogenized coefficient is of about 10\%. In contrast, we plot (see the red curve) the function $\overline{k} \mapsto J(\overline{k})$ where $J$ is again defined by~\eqref{eq:def_J}, but where $(\overline{u},u_\eps)$ is now the solution to~\eqref{eq:pb_min_H} (or equivalently~\eqref{eq:arlequin}) with $W_H$ replaced by $W_H^{\rm enrich}$. We see that this function attains its minimum for a value of $\overline{k} \approx 1.42$ very close to $k^\star$: the error on the evaluation of the homogenized coefficient is reduced to 0.4\%. At no practical additional cost (since the only additional computational cost is the offline precomputation of the solution to~\eqref{eq:def_psi0}, and the online addition of {\em one} dimension in the discretization space for $\psi$), we hence obtain a much more accurate approximation of the homogenized coefficient when considering $W_H^{\rm enrich}$. In all what follows, and without writing it explicitly, we use the space $W_H^{\rm enrich}$ rather than $W_H$.

\medskip

\subsection{On the choice of the initial guess for the optimization loop} \label{sec:IG}

Solving the optimization problem~\eqref{eq:optim_J} amounts to solving several times the equations~\eqref{eq:arlequin}. Each resolution of~\eqref{eq:arlequin} is expensive because the unknown function $u_\eps$ is discretized on the fine mesh ${\cal T}_h$ (with a mesh size $h \ll \eps$, say $h \approx \eps/10$). The number of degrees of freedom involved in~\eqref{eq:arlequin} is thus large. It is therefore interesting to start the optimization procedure~\eqref{eq:optim_J} from a good initial guess, i.e. an initial guess as close as possible to the minimizer. 

Consider for instance (this is just an example and our arguments carry over to other cases, as pointed out in Remark~\ref{rem:penal} below) a discretization strategy where the Dirichlet boundary condition $\overline{u}(x) = x_1$ on $\Gamma$ is taken into account by penalization. This is the case in the software we use, namely FreeFem++~\cite{hecht}. In practice, instead of~\eqref{eq:arlequin}, we solve for $\overline{u} \in V_H$, $u_\eps \in V_h$ and $\psi \in W_H$ such that
\begin{equation} \label{eq:arlequin_pena}
  \begin{cases}
  \dps \forall \overline{v} \in V_H, \quad \overline{A}(\overline{u},\overline{v}) + \frac{1}{\eta} \int_{\Gamma} (\overline{u} - x_1) \, \overline{v} + C(\overline{v},\psi) = 0,
  \\
  \forall v_\eps \in V_h, \quad A_\eps(u_\eps,v_\eps) - C(v_\eps,\psi) = 0,
  \\
  \forall \phi \in W_H, \quad C(\overline{u}-u_\eps,\phi) = 0,
  \end{cases}
\end{equation}
where $\eta$ is a small positive parameter. Note that, when $\eta \to 0$, the solution to~\eqref{eq:arlequin_pena} converges to that of~\eqref{eq:arlequin}.

The equation~\eqref{eq:arlequin_pena} takes the form
\begin{equation} \label{eq:arlequin_pena_matrix}
\left( \begin{array}{ccc}
\overline{\cal A} + \overline{\cal P} & 0 & \overline{\cal C} \\
0 & {\cal A}_\eps & -{\cal C}_\eps \\
\overline{\cal C}^T & -{\cal C}_\eps^T & 0 
\end{array} \right)
\
\left( \begin{array}{c}
\overline{U} \\ U_\eps \\ \Psi 
\end{array} \right)
=
\left( \begin{array}{c}
\overline{F} \\ 0 \\ 0
\end{array} \right),
\end{equation}
where $\overline{\cal A}$ (resp. ${\cal A}_\eps$) is the matrix representing the bilinear form $\overline{A}$ (resp. $A_\eps$) in the finite dimensional space $V_H \times V_H$ (resp. $V_h \times V_h$). The matrix $\overline{\cal C}$ (resp. ${\cal C}_\eps$) represents the bilinear form $C$ in the finite dimensional space $V_H \times W_H$ (resp. $V_h \times W_H$). The matrix $\overline{\cal P}$ represents the bilinear form $\dps (\overline{u},\overline{v}) \mapsto \eta^{-1} \int_{\Gamma} \overline{u} \ \overline{v}$. The vector $\overline{F}$ on the right-hand side represents the linear form $\dps \overline{v} \mapsto \eta^{-1} \int_{\Gamma} x_1 \, \overline{v}$. The unknown functions $\overline{u} \in V_H$, $u_\eps \in V_h$ and $\psi \in W_H$ are represented by the vectors $\overline{U}$, $U_\eps$ and $\Psi$, respectively.

Since $\overline{\cal A}$ depends linearly on the scalar $\overline{k}$ while $\overline{F}$, $\overline{\cal P}$, ${\cal A}_\eps$, $\overline{\cal C}$ and ${\cal C}_\eps$ are independent of $\overline{k}$, the above linear system~\eqref{eq:arlequin_pena_matrix} can be recast as
\begin{equation} \label{eq:sys_lin}
\left( \overline{k} Z_1 + Z_2 \right) U(\overline{k}) = F_0,
\end{equation}
where the vector $F_0$ and the matrices $Z_1$ and $Z_2$ do not depend on $\overline{k}$ (the vector $U$ of course collects the vectors $\overline{U}$, $U_\eps$ and $\Psi$).

\begin{remark} \label{rem:penal}
To take into account the Dirichlet boundary condition $\overline{u}(x) = \overline{u}_{\rm ref}(x) = x_1$ on $\Gamma$, another possibility (alternative to penalization) is to change the unknown function and to work with $\widetilde{u} = \overline{u} - \overline{u}_{\rm ref}$. Instead of~\eqref{eq:arlequin}, we then solve for $\widetilde{u} \in V_H^0$, $u_\eps \in V_h$ and $\psi \in W_H$ such that
$$
\begin{cases}
  \dps \forall \overline{v} \in V_H^0, \quad \overline{A}(\widetilde{u},\overline{v}) + C(\overline{v},\psi) = -\overline{A}(\overline{u}_{\rm ref},\overline{v}),
  \\
  \forall v_\eps \in V_h, \quad A_\eps(u_\eps,v_\eps) - C(v_\eps,\psi) = 0,
  \\
  \forall \phi \in W_H, \quad C(\widetilde{u}-u_\eps,\phi) = -C(\overline{u}_{\rm ref},\phi).
\end{cases}
$$
This problem yields a linear system where the matrix again depends on $\overline{k}$ in an affine way. In contrast with the penalization formulation~\eqref{eq:arlequin_pena_matrix}, the right-hand side now also depends on $\overline{k}$ in an affine way. The argument below can be adapted to that case. 
\end{remark}
  
Suppose now that we have at our disposal the solution $U(\overline{k}_0)$ to~\eqref{eq:sys_lin} for some $\overline{k}_0$. We thus write
$$
\left( \overline{k} Z_1 + Z_2 \right) U(\overline{k}) = \left( \overline{k}_0 Z_1 + Z_2 \right) U(\overline{k}_0),
$$
hence
\begin{equation} \label{eq:moto2}
M(\overline{k}) \, U(\overline{k}) = U(\overline{k}_0),
\end{equation}
where the matrix $M(\overline{k})$ depends on $\overline{k}$ in an affine manner:
\begin{align} 
  M(\overline{k}) &= \overline{k} \left( \overline{k}_0 Z_1 + Z_2 \right)^{-1} Z_1 + \left( \overline{k}_0 Z_1 + Z_2 \right)^{-1} Z_2
  \nonumber
  \\
  &= (\overline{k} - \overline{k}_0) \left( \overline{k}_0 Z_1 + Z_2 \right)^{-1} Z_1 + \Id.
  \label{eq:def_M}
\end{align}
Here $\Id$ is an identity matrix of the size  $N \times N$ with  $N = \text{dim } V_H + \text{dim } V_h +\text{dim } W_H$. 

Suppose temporarily that we wish to find $\overline{k}$ such that $U(\overline{k}) = U_{\rm target}$ for some given $U_{\rm target}$. Finding such $\overline{k}$ is not easy because $U(\overline{k})$ depends on $\overline{k}$ in a complex manner. However, this problem is equivalent to finding $\overline{k}$ such that $M(\overline{k}) \, U(\overline{k}) = M(\overline{k}) \, U_{\rm target}$, which is itself equivalent to $U(\overline{k}_0) = M(\overline{k}) \, U_{\rm target}$. We recast the problem as finding $\overline{k}$ that minimizes $\| U(\overline{k}_0) - M(\overline{k}) \, U_{\rm target} \|^2$. This formulation is much easier to solve, because of the affine formula~\eqref{eq:def_M}. Setting $Z_0 = \left( \overline{k}_0 Z_1 + Z_2 \right)^{-1}$, we obtain
\begin{equation} \label{eq:IG_k1}
\overline{k} = \overline{k}_0 + \frac{\big( U(\overline{k}_0) - U_{\rm target} \big) \cdot Z_0 Z_1 U_{\rm target}}{\| Z_0 Z_1 U_{\rm target} \|^2}.
\end{equation}

In practice, the quantity $J$ that we minimize only depends on $\overline{u}$, and not on $(\overline{u},u_\eps,\psi)$. We thus cannot think of our problem as finding some $\overline{k}$ such that $U(\overline{k}) = U_{\rm target}$ and we cannot directly use~\eqref{eq:IG_k1}, since we do not know $U_{\rm target}$, but only its first component (namely the function $\overline{u}_{\rm ref}$). Therefore, we cannot expect the above formula~\eqref{eq:IG_k1} to yield the minimizer of $J$. It however can serve as a good starting point to build an adequate initial guess for the optimization problem~\eqref{eq:optim_J}. For this reason, we consider as initial guess the formula
\begin{equation} \label{eq:IG_k2}
\overline{k} = \overline{k}_0 + \frac{\pi \big( U(\overline{k}_0) - U_{\rm target} \big) \cdot Z_0 Z_1 \pi \big( U_{\rm target} \big)}{\| Z_0 Z_1 \pi \big( U_{\rm target} \big) \|^2},
\end{equation}
with, for any vectors $\overline{U}$, $U_\eps$ and $\Psi$,
$$
\pi \left( \begin{array}{c}
\overline{U} \\ U_\eps \\ \Psi 
\end{array} \right)
=
\left( \begin{array}{c}
\overline{U} \\ 0 \\ 0
\end{array} \right).
$$
As we now show, this choice~\eqref{eq:IG_k2} of initial guess allows to significantly speed up the optimization algorithm solving~\eqref{eq:optim_J}. We also note that~\eqref{eq:IG_k2} can be extended to the case of a matrix-valued coefficient $\overline{k}$, as detailed in Appendix~\ref{sec:2Dformula}.

\medskip

\begin{figure}[htbp]
\centering
\includegraphics[width=0.6\textwidth]{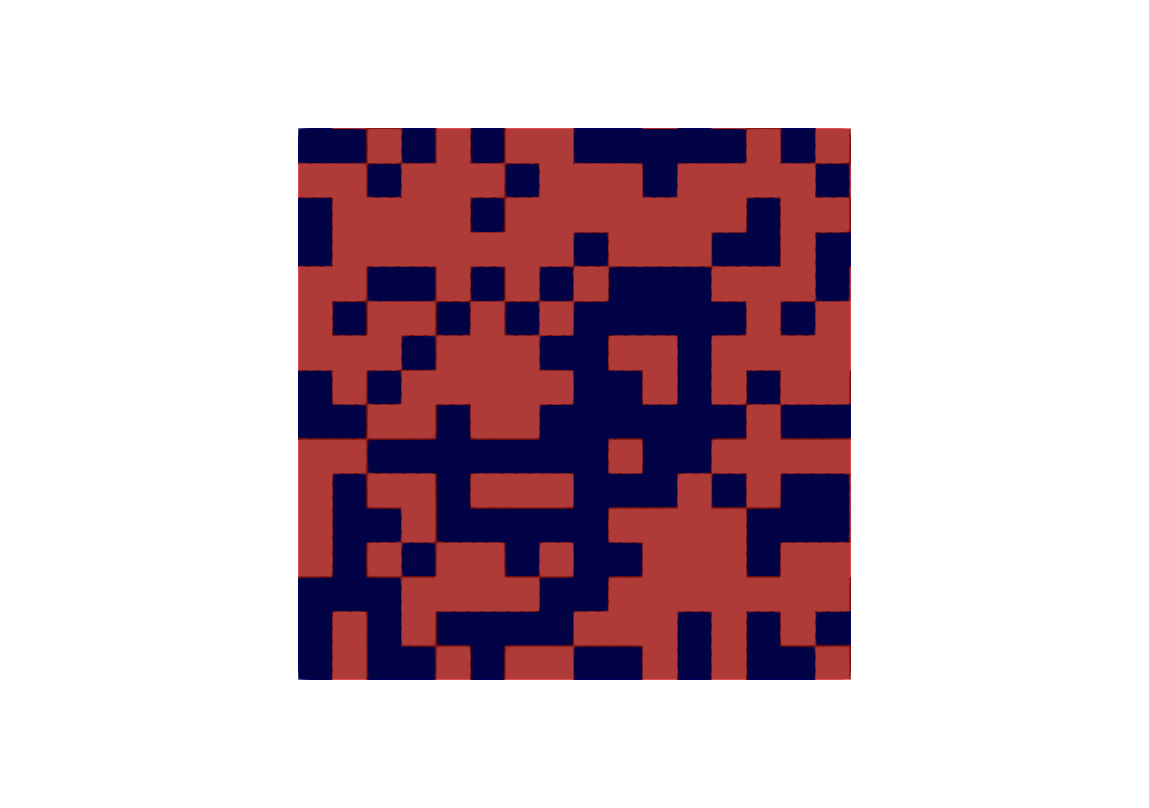}
\caption{Representation of the function $x \mapsto k_\eps(x)$. Each square is of size $\eps \times \eps$. On the red squares, $k_\eps(x) = 1$. On the blue squares, $k_\eps(x) = 2$. \label{fig:board}}
\bigskip
\end{figure}

To illustrate the efficiency of~\eqref{eq:IG_k2}, we consider the following numerical example. On $D_c \cup D_f$, we consider the heterogeneous coefficient shown on Figure~\ref{fig:board}. We set $\eps = 1/16$, and use a fine mesh of size $h = 1/128$. The coarse mesh size is $H = 1$. 

\medskip

The optimization problem~\eqref{eq:optim_J} is solved using the Newton algorithm which is, as is well-known, all the more so efficient as it starts from an initial guess (IG) close to the solution. We therefore expect that using the initial guess~\eqref{eq:IG_k2} significantly reduces the number of iterations of the Newton algorithm, and thus its cost. This is indeed the case. If we start with the initial guess $\overline{k}^{\rm IG} = 0.8$, the algorithm needs 9 iterations to reach a prescribed tolerance (relative variation of $J$ smaller than $10^{-2}$). In contrast, if we start from the initial guess~\eqref{eq:IG_k2} computed with $\overline{k}_0 = 0.8$ (namely $\overline{k}^{\rm IG} \approx 1.419$), we only need 2 iterations to reach the same tolerance.

For the sake of completeness, we have also considered solving the optimization problem~\eqref{eq:optim_J} using the Nelder-Mead algorithm, which is a zero-order algorithm that does not need to compute the gradient of $J$ with respect to $\overline{k}$. The Nelder-Mead algorithm is very efficient at globally exploring the functional to optimize. In contrast, when initialized close to a local optimum, it is less efficient than first or second-order algorithms (including the Newton algorithm). In particular, it does not necessarily stay close to that local optimum. It is interesting to observe that using the initial guess~\eqref{eq:IG_k2} is also beneficial when using the Nelder-Mead algorithm (although the gain is of course less significant than with the Newton algorithm, a fact which could be expected given the fact that the {\em purpose} of the two algorithms are different, as we have recalled above): instead of 18 iterations when starting from $\overline{k}^{\rm IG} = 0.8$, the Nelder-Mead algorithm needs 16 iterations when starting from the initial guess $\overline{k}^{\rm IG} \approx 1.419$ given by~\eqref{eq:IG_k2}.

\subsection{Solution procedure for the random case} \label{sec:random_nous}

We now return to the random case considered in Section~\ref{sec:random_cottereau}, where we have pointed out that the major difficulty is to solve the huge system~\eqref{eq:arlequin_sto}. We now explain how we believe that this problem can be expressed in a way amenable to computations at a limited cost.

The second line of~\eqref{eq:arlequin_sto} reads as: for any $1 \leq m \leq M$, find $u_\eps^m \in V_h$ such that
\begin{equation} \label{eq:neumann}
\forall v_\eps^m \in V_h, \quad A_\eps^m(u_\eps^m,v_\eps^m) = C(v_\eps^m,\psi).
\end{equation}
This is a Neumann problem on $u_\eps^m$, and the above equation hence determines (once $\psi$ is known) the function $u_\eps^m$ up to the addition of a constant. We thus write
$$
  u_\eps^m = \beta^m + {\cal L}_\eps^m \psi,
$$
where $\beta^m$ is a constant and ${\cal L}_\eps^m \psi$ is the solution of vanishing average of~\eqref{eq:neumann}. Note that the problems~\eqref{eq:neumann}, for $1 \leq m \leq M$, are independent from each other: they each amount to solving for a {\em single} fine-scale function, and not $M$ such functions as in~\eqref{eq:arlequin_sto}. 

We see that~\eqref{eq:arlequin_sto} actually implies that $\overline{u} \in V_H^{\rm Dir BC}$ and $\psi \in W_H$ satisfy the following problem:
\begin{equation} \label{eq:arlequin_sto_compact}
  \begin{cases}
  \dps \forall \overline{v} \in V_H^0, \quad \overline{A}(\overline{u},\overline{v}) + C(\overline{v},\psi) = 0,
  \\
  \dps \forall \phi \in W_H, \quad C\left(\overline{u} - \frac{1}{|D_c|} \int_{D_c} \overline{u} - \frac{1}{M} \sum_{m=1}^M {\cal L}_\eps^m \psi, \phi \right) = 0,
  \end{cases}
\end{equation}
where the number of degrees of freedom is limited (and independent of $M$ and $h$). The problem~\eqref{eq:arlequin_sto_compact} yields a linear system of the form
\begin{equation} \label{eq:condens}
\left( \begin{array}{cc}
\overline{\cal A} & \overline{\cal C} \\
\widetilde{\cal C}^T & \quad \dps - \frac{1}{M} \sum_{m=1}^M {\cal C}_\eps^T \left( {\cal A}_\eps^m \right)^{-1} {\cal C}_\eps
\end{array} \right)
\
\left( \begin{array}{c}
\overline{U} \\ \Psi 
\end{array} \right)
=
\left( \begin{array}{c}
\overline{F} \\ 0
\end{array} \right),
\end{equation}
where $\overline{\cal A}$ is the matrix representing the bilinear form $\overline{A}$ in the finite dimensional space $V_H^0 \times V_H^0$, $\overline{\cal C}$ is the matrix representing the bilinear form $C$ in $V_H^0 \times W_H$, $\widetilde{\cal C}$ is the matrix representing the bilinear form $\dps (\overline{v},\phi) \mapsto C\left(\overline{v} - \frac{1}{|D_c|} \int_{D_c} \overline{v}, \phi \right)$ in $V_H^0 \times W_H$, ${\cal C}_\eps$ is the matrix representing the bilinear form $C$ in $V_h \times W_H$ (as in~\eqref{eq:arlequin_pena_matrix}) and $\dps \left( {\cal A}^m_\eps \right)^{-1}$ is the inverse of the matrix ${\cal A}^m_\eps$ representing the bilinear form $A^m_\eps$ in $V_h \times V_h$ (as in~\eqref{eq:arlequin_pena_matrix}). More precisely, $\dps \left( {\cal A}^m_\eps \right)^{-1}$ is a specific inverse of the matrix ${\cal A}^m_\eps$ in the sense that, for any vector $F$ representing a function of vanishing mean, $U = \left( {\cal A}^m_\eps \right)^{-1} F$ means that ${\cal A}^m_\eps U = F$ and that $U$ represents a function of vanishing mean (otherwise stated, $\left( {\cal A}^m_\eps \right)^{-1}$ is the matrix representing the operator ${\cal L}_\eps^m$ defined below~\eqref{eq:neumann}). The number of unknowns in~\eqref{eq:condens} is small: it is simply equal to $\text{dim } V_H + \text{dim } W_H$, and thus only depends on the coarse mesh size $H$ and the dimension $d$.

\medskip

We now turn to a numerical example that illustrates this procedure. In order to show the versality of the approach, we consider here the case when the tentative constant coefficient $\overline{k}$ upon which we optimize (in order to identify $k^\star$) is matrix-valued, in contrast to the case considered so far where we have assumed $\overline{k}$ to be a scalar. The optimization approach described in Section~\ref{sec:optim} can be extended to that case. Lemma~\ref{lem:consistence} is then modified as follows: (i) if $\overline{k} \, e_1 = k^\star \, e_1$, then the solution to~\eqref{eq:pb_min_star} is $\overline{u}(x) = x_1$ in $D \cup D_c$ and $u_0(x) = x_1$ in $D_c \cup D_f$; (ii) conversely, if $(\overline{u},u_0)$ is a solution to~\eqref{eq:pb_min_star} with $\overline{u}(x) = x_1$ in $D \cup D_c$, then $u_0(x) = x_1$ in $D_c \cup D_f$ and $\overline{k} \, e_1 = k^\star \, e_1$. The optimization problem~\eqref{eq:optim_J}--\eqref{eq:def_J} can hence only be used to recover $k^\star \, e_1$, and not the complete matrix $k^\star$.

When $k^\star$ is matrix-valued, we thus need to consider two optimization problems to recover $k^\star$: first Problem~\eqref{eq:optim_J}--\eqref{eq:def_J}, and second a similar one where $\dps J(\overline{k}) = \int_{D \cup D_c} | \nabla \overline{u} - e_2 |^2$, where $(\overline{u},u_0)$ is the solution to~\eqref{eq:pb_min_star} where the boundary condition $\overline{u}(x) = x_1$ on $\Gamma$ is of course replaced by $\overline{u}(x) = x_2$ on the part of the exterior boundary of $D$ made by the two horizontal lines (that is $\partial D \setminus (\Gamma \cup \partial D_c)$, see Figure~\ref{fig:decompo_left}).

\medskip

We consider the random function
\begin{equation} \label{eq:RCB}
k(x,\omega) = \sum_{j \in \ZZ^2} k_j(\omega) \ 1_{j + Q}(x),
\end{equation}
where $Q=(0,1)^2$ is the unit square, and where $k_j$ are i.i.d. random variables. We next set $k_\eps(x,\omega) = k(x/\eps,\omega)$. Figure~\ref{fig:board} shows a particular realization of the function $x \mapsto k_\eps(x,\omega)$. We set $\eps = 1/32$ and assume that $k_j(\omega)$ is equal to either 1 or 64 with a probability $1/2$. The fine mesh size is $h = 1/256$, and the coarse mesh size is $H = 1/2$. We use $M=96$ realizations to discretize the expectation. 

In that case, the homogenized coefficient is known (see e.g.~\cite{le2012}): it is proportional to the identity matrix and is equal to $k^\star = 8 \ \Id$. Using as boundary conditions $\overline{u}(x) = x_1$ on $\Gamma$ and optimizing $J$ with respect to the coefficient $\overline{k} \in \RR^{2 \times 2}$, we are only in position to identify $\overline{k}_{11}$ and $\overline{k}_{21}$, as explained above. In that case, the minimum of $J$ is attained when $\overline{k}_{11} = 7.9842$ and $\overline{k}_{21} = 0.0123$. We thus obtain a very accurate approximation of the actual homogenized coefficient, for a limited computational cost.

\subsection{Selection of the random configurations} \label{sec:sqs}

Still considering the random setting of Section~\ref{sec:random_nous}, we now show how variance reduction type methods can be implemented within the approach described here, in order to reduce the statistical noise in the approximation of $k^\star$.

In practice, as shown above, only a finite number $M$ of realizations can be considered to approximate the expectations. The expectations are thus replaced by empirical averages, which are themselves random quantities. The output $\overline{k}^{\rm opt}(\eps)$ of the optimization algorithm is therefore random. This output of course accurately approximates the deterministic limit when $M$ is taken extremely large, but the computational cost of the whole procedure then increases prohibitively. It is however possible to reduce the statistical noise on the output, at essentially no extra computational cost (and thus, for a fixed accuracy of the final result, reduce the computational cost), by using the following selection procedure.

Instead of considering $M$ arbitrary realizations of $k_\eps(x,\omega)$ in~\eqref{eq:arlequin_sto}, we suggest to consider $M$ well-chosen realizations, in the spirit of our previous work~\cite{sqs}. These selected realizations (which we call hereafter the Special Quasirandom Structures (SQS)) are expected, for a finite value of $\eps$, to better represent the asymptotic limit $\eps \to 0$ than a set of $M$ generic realizations. In a slightly different context, significant computational gains have been obtained in~\cite{sqs} when using this selection approach.

More precisely, we suggest to consider realizations $\omega$ such that the difference
\begin{equation} \label{eq:sqs_diff}
\frac{1}{|D_c \cup D_f|} \int_{D_c \cup D_f} k_\eps(x,\omega) \, dx - \frac{1}{|D_c \cup D_f|} \int_{D_c \cup D_f} \EE[k_\eps(x,\cdot)] \, dx,
\end{equation}
which is easy to evaluate, is as small as possible (see~\cite[Sec.~2.3.2]{sqs}). 

\begin{remark}
For the case of the random checkerboard considered in Section~\ref{sec:random_nous} (see~\eqref{eq:RCB}), where $k_j(\omega)$ is either equal to $k^0$ or $k^1$ with equal probability 1/2, the difference~\eqref{eq:sqs_diff} is directly related to the difference between the ratio of cells $k^0$ vs $k^1$ and 1, which is the asymptotic ratio in the limit of infinitely many cells. Enforcing that~\eqref{eq:sqs_diff} is small thus amounts to enforcing that the ratio of cells $k^0$ vs $k^1$ is as close as possible to one (its asymptotic limit when $\eps \to 0$). 

We also note that, in the specific case of the random checkerboard, it is actually possible to draw realizations such that~\eqref{eq:sqs_diff} vanishes (that is, realizations with an equal number of cells $k^0$ and $k^1$).
\end{remark}

The selection algorithm thus proceeds as follows. 

\begin{algo}[{\bf Selection algorithm}] \label{algo:selection}
Fix a number ${\cal M} \gg M$ of trials.
\begin{enumerate}
\item Compute the second term of~\eqref{eq:sqs_diff}.
\item For $m = 1, \dots, {\cal M}$,  
\begin{enumerate}
\item Generate a random environment $\omega_m$.
\item Compute the first term of~\eqref{eq:sqs_diff}.
\item Compute the error ${\tt error}_m$ between the first and the second terms of~\eqref{eq:sqs_diff}.
\end{enumerate}
\item Sort the random environments $(\omega_m)_{1 \leq m \leq {\cal M}}$ according to ${\tt error}_m$. Keep the best $M$ realizations (and use these in~\eqref{eq:arlequin_sto} or equivalently~\eqref{eq:arlequin_sto_compact}), and reject the others.
\end{enumerate}
\end{algo}

We note that, in full generality, the cost of Monte Carlo approaches is usually dominated by the cost of draws, and therefore selection algorithms are targeted to reject as few draws as possible. In the present context, where boundary value problems such as~\eqref{eq:neumann} are to be solved repeatedly (because of the many realizations, and because of the optimization context), the cost of draws for the environment is negligible in comparison to the cost of the solution procedure for such boundary value problems. Likewise, evaluating both terms of~\eqref{eq:sqs_diff} is not expensive. Therefore, the purpose of the selection mechanism is to limit the number of boundary value problems to be solved, even though this comes at the (tiny) price of rejecting many environments. This explains why we employ a simplistic rejection procedure for the selection, while in other situations of Monte Carlo samplings, one would invest in a more elaborate selection procedure.

\medskip

We now consider a numerical example to illustrate the above approach. Let $I \in \NN^\star$ be given. We are going to compare the following two computational strategies, which essentially share the same cost:
\begin{enumerate}
\item Reference Monte Carlo approach:
  \begin{itemize}
  \item for any $1 \leq i \leq I$, perform the following computation:
    \begin{itemize}
    \item draw $M$ realizations of the function $k_\eps(x,\omega)$ that are identically distributed, independent from each other and independent of those drawn for the previous values of $i$;
    \item compute the optimal constant coefficient $\overline{k}^{{\rm opt},i}_{\rm MC}$ by considering~\eqref{eq:arlequin_sto} with these $M$ heterogeneous functions;
    \end{itemize}
  \item compute the variance $V_{\rm MC}$ of the outputs $\{ \overline{k}^{{\rm opt},i}_{\rm MC} \}_{1 \leq i \leq I}$.
  \end{itemize}
\item Selection approach:
  \begin{itemize}
  \item for any $1 \leq i \leq I$, perform the following computation:
    \begin{itemize}
    \item draw ${\cal M}$ realizations of the function $k_\eps(x,\omega)$ that are identically distributed, independent from each other and independent of those drawn for the previous values of $i$;
    \item select the best $M$ realizations among these ${\cal M}$ realizations according to Algorithm~\ref{algo:selection};
    \item compute the optimal constant coefficient $\overline{k}^{{\rm opt},i}_{\rm SQS}$ by considering~\eqref{eq:arlequin_sto} with these $M$ best heterogeneous functions;
    \end{itemize}
  \item compute the variance $V_{\rm SQS}$ of the outputs $\{ \overline{k}^{{\rm opt},i}_{\rm SQS} \}_{1 \leq i \leq I}$.
  \end{itemize}
\end{enumerate}
In practice, we again consider the random heterogeneous function $k_\eps(x,\omega) = k(x/\eps,\omega)$ where $k$ is defined by~\eqref{eq:RCB} and $\eps = 1/16$. We assume that $k_j(\omega)$ is equal to either 1 or some $\beta > 1$ with a probability $1/2$. The parameter $\beta$ hence represents the contrast of the random checkerboard. We compute the variance from $I=100$ different optimal constant coefficients $\overline{k}^{{\rm opt},i}$, $1 \leq i \leq I$. The computation of each $\overline{k}^{{\rm opt},i}$ is performed using the Arlequin approach (see~\eqref{eq:arlequin_sto_compact}) using $M=16$ random environments and the mesh sizes $h=1/128$ and $H=1$. 

\medskip

Our results are collected in Table~\ref{tab:variance}. When $\beta = 1.5$, the selection approach provides results the variance of which is (almost) 5 times as small as when using the Monte Carlo approach. This is a significant computational gain. Even for a comparably large value of the contrast, namely $\beta = 64$, we observe that the selection procedure decreases the variance by a factor of 2. As always, the gain in variance decreases when the contrast increases. This observation is consistent with observations already made in~\cite{sqs} where the variance reduction approach is applied for the approximation of the homogenized coefficients. However, in that work~\cite{sqs}, the selection procedure leads to significantly larger computational gains than here. A possible explanation of this effect is that the Arlequin approach may be more stable with respect to the statistical noise than the classical homogenization procedure, since it avoids a direct computation of the corrector functions.

\bigskip

\begin{table}[htbp]
\centering{
\begin{tabular}{c@{\hskip 3mm}c@{\hskip 5mm}c@{\hskip 5mm}c} 
\hline
& & & \\
Contrast $\beta$ & $V_{\rm MC}$ & $V_{\rm SQS}$ & $V_{\rm MC}/ V_{\rm SQS}$ \\[3mm]
\hline
\hline
\midrule
1.5 & $9.863 \times 10^{-7}$ & $2.242 \times 10^{-7}$ & 4.40 \\
9   & $2.693 \times 10^{-4}$ & $1.192 \times 10^{-4}$ & 2.26 \\
64  & $2.152 \times 10^{-2}$ & $1.125 \times 10^{-2}$ & 1.91 \\
\midrule
\hline \\
\end{tabular}}
\caption{Variances of the result without and with the selection procedure, for different values of the contrast in $k_\eps$. \label{tab:variance}}
\bigskip
\end{table}

\subsection{Computation of the corrector} \label{sec:corrector}

All the above sections focus on approximating the homogenized coefficient $k^\star$. However, in the homogenization context, constructing an approximation of the corrector function is equally interesting, because it leads to an accurate approximation of the solution (in the natural $H^1$ norm, i.e. the energy norm) of the oscillatory problem. We thus wish to extend our approach in that direction. Besides its own interest, such an extension also allows for a fair comparison between our approach here and the classical homogenization approach, which yields both the corrector function and the homogenized matrix. 

To illustrate our strategy, we first consider a problem where only one corrector function is to be determined, say $w_1$ (more general cases are considered below). The lamellar case, where $k_\eps(x_1,x_2) = k_{\rm per}(x_1/\eps)$ for some periodic function $k_{\rm per}$, falls within that framework. In addition, $w_1$ only depends on $x_1$. For a standard heterogeneous problem of the form $-\dive (k_\eps \nabla \widetilde{u}_\eps) = f$ in some domain, it is well-known that the oscillatory solution can be approximated by a two-scale expansion written in terms of the corrector function $w_1$ and the homogenized limit $u^\star$ (see e.g.~\cite{le2012}): 
\begin{equation} \label{eq:tse_1correc}
\widetilde{u}_\eps(x_1,x_2) \approx u^\star(x_1, x_2) + \eps w_1(x_1/\eps) \, \partial_{x_1} u^\star(x_1, x_2).
\end{equation}
We then replace $(\widetilde{u}_\eps, u^\star)$ with the solution $(u_\eps, \overline{u})$ to the Arlequin approach (namely, the solution to~\eqref{eq:pb_min_H}) when $\overline{k}$ has reached its optimal value (which is close to $k^\star$).
Following an idea already used in~\cite{cras_kun_li}, we thus expect to recover $w_1$ by writing
\begin{equation}\label{eq:corrector_approximated}
w_1'\left(\frac{x_1}{\eps}\right) \approx \widetilde{W}_\eps(x_1,x_2) := \frac{\partial_{x_1} u_\eps(x_1,x_2) - \partial_{x_1} \overline{u}(x_1,x_2)}{\partial_{x_1} \overline{u}(x_1,x_2)}.
\end{equation}
The ratio~\eqref{eq:corrector_approximated} can be evaluated in $D_c$ where the two models co-exist (note that $\overline{u}(x_1,x_2)$ may also be replaced by $x_1$, which is its asymptotic value when $\eps \to 0$ for the optimal choice of $\overline{k}$), or in $D_f$ (in that case, we replace $\overline{u}(x_1,x_2)$ by $x_1$ in the above expression). In passing, we note that $k^\star$ is not scalar-valued for this lamellar case. The enrichment of the Lagrange multiplier space by the solution $\psi_0$ to~\eqref{eq:def_psi0} is however still sufficient, since $k^\star e_1$ is colinear to $e_1$ in this case.

\medskip

We start from the following simple numerical illustration for the lamellar case. We set
$$
k_\eps(x_1,x_2) = k_{\rm per}(x_1/\eps) =  \left(1+\sin^2(2 \pi x_1/\eps) \right)^{-1},
$$
in which case the corrector function is analytically known:
$$
w_1'\left(\frac{x_1}{\eps}\right) = \frac{2}{3 \, k_{\rm per}(x_1/\eps)} - 1.
$$
In Table~\ref{tab:correctors}, we compare two approximations of $w_1'(x_1/\eps)$ in the case $\eps = 1/16$. The first one (see the left part of~\eqref{eq:def_error}) is provided by $\widetilde{W}_\eps(x_1,x_2)$ defined by~\eqref{eq:corrector_approximated} in $D_c$, where $(\overline{u},u_\eps)$ is the solution to~\eqref{eq:pb_min_H} for some discretization parameters $h$ and $H$. The second one (see the right part of~\eqref{eq:def_error}) is provided by a numerical computation of the corrector function using the classical homogenization approach (see for instance~\cite{le2012}). We hence solve the corrector equation at the scale $\eps$, using a fine mesh of size $h$, and therefore compute some $w_1^h$. To compare these two approximations, we define the errors
\begin{equation} \label{eq:def_error}
  \widetilde{e} = \left\| \widetilde{W}_\eps - w_1'\left(\frac{\cdot}{\eps}\right) \right\|_{L^2(D_c)}
  \quad \text{and} \quad
  e = \left\| \partial_1 w_1^h\left(\frac{\cdot}{\eps}\right) - w_1'\left(\frac{\cdot}{\eps}\right) \right\|_{L^2(D_c)}.
\end{equation}
We observe in Table~\ref{tab:correctors} that the error $\widetilde{e}$ decreases linearly in terms of $H$, and that our procedure yields an approximation $\widetilde{W}_\eps$ which is as accurate as the classical homogenization approach when $H$ is small enough (say here $H \leq 1/8$): in that case, $\widetilde{e} \approx e$.

Of course, for a general case, we do not have access to the exact corrector, and we therefore cannot compute $\widetilde{e}$. This is the reason why, on Table~\ref{tab:correctors}, we also show the error
\begin{equation} \label{eq:def_error_bis}
  \widetilde{e}_{\rm approx} = \left\| \widetilde{W}_\eps - \partial_1 w_1^h\left(\frac{\cdot}{\eps}\right) \right\|_{L^2(D_c)}.
\end{equation}
In the range of parameters that we consider, we observe that the error~\eqref{eq:def_error_bis} is very close to the actual error $\widetilde{e}$ of our approach, and therefore an accurate estimate of it. In the numerical test described below, for which we do not have access to the exact correctors, we will therefore rely on~\eqref{eq:def_error_bis} to assess the accuracy of our approximation. 

\bigskip

\begin{table}[htbp]
\centering{
\begin{tabular}{c@{\hskip 5mm}c@{\hskip 5mm}c@{\hskip 5mm}c@{\hskip 5mm}c} 
\hline
& & \\
$H$ & $\widetilde{e}$ & $\widetilde{e}_{\rm approx}$ & $\dps \frac{\widetilde{e}}{\left\| w_1'(\cdot/\eps) \right\|_{L^2(D_c)}}$ & $\dps \frac{\widetilde{e}_{\rm approx}}{\left\| \partial_1 w_1^h(\cdot/\eps) \right\|_{L^2(D_c)}}$ \\[3mm]
\hline
\hline
\midrule
1/16 & 0.00647 & 0.00432 & 0.0345 & 0.0229 \\
1/8  & 0.00839 & 0.00661 & 0.0448 & 0.0364 \\
1/4  & 0.0136  & 0.0125  & 0.0725 & 0.0675 \\
1/2  & 0.0251  & 0.0244  & 0.133  & 0.130  \\
1    & 0.0454  & 0.0447  & 0.240  & 0.237  \\ 
\midrule
\hline \\
\end{tabular}}
\caption{Values of the absolute errors~\eqref{eq:def_error} and~\eqref{eq:def_error_bis} (and corresponding relative errors) for different sizes $H$ of the coarse mesh ($\eps = 1/16$ and $h = 1/128$). The absolute error when using the classical homogenization approach is $e = 0.00509$ (which corresponds to a relative error $e/\left\| w_1'(\cdot/\eps) \right\|_{L^2(D_c)} = 0.0272$). \label{tab:correctors}}
\bigskip
\end{table}

The above results suggest the following strategy. We first determine an accurate approximation $\overline{k}^{\rm opt}$ of $k^\star$ by considering the optimization problem~\eqref{eq:optim_J}--\eqref{eq:def_J}, where the Arlequin equations~\eqref{eq:pb_min_H} are solved using a fine mesh of size $h \ll \eps$ adapted to the oscillations of $k_\eps$ and a coarse mesh of size $H_{\rm coarse}$ that we choose to be relatively large (say $H_{\rm coarse} = 1$). We have indeed observed in Section~\ref{sec:enrich} that it is possible to obtain an accurate approximation of $k^\star$ when using the coarse value $H=1$ (recall that $D \cup D_c \cup D_f = (-4,4)^2$, so this value of $H$ corresponds to only 8 elements per direction in the computational domain). Of course, since determining an approximation of $k^\star$ is an iterative problem, it is very beneficial to use a value of $H$ as large as possible. Next, in a second stage, we solve one more time the equations~\eqref{eq:pb_min_H}, with the obtained value $\overline{k}^{\rm opt}$ of the constant coefficient, but this time using a coarse mesh of size $H_{\rm fine} \ll H_{\rm coarse}$ (say $H_{\rm fine} = 1/8$), and compute~\eqref{eq:corrector_approximated} using the latter solution. 

\medskip

\begin{figure}[htbp]
\centering
\includegraphics[width=0.48\textwidth]{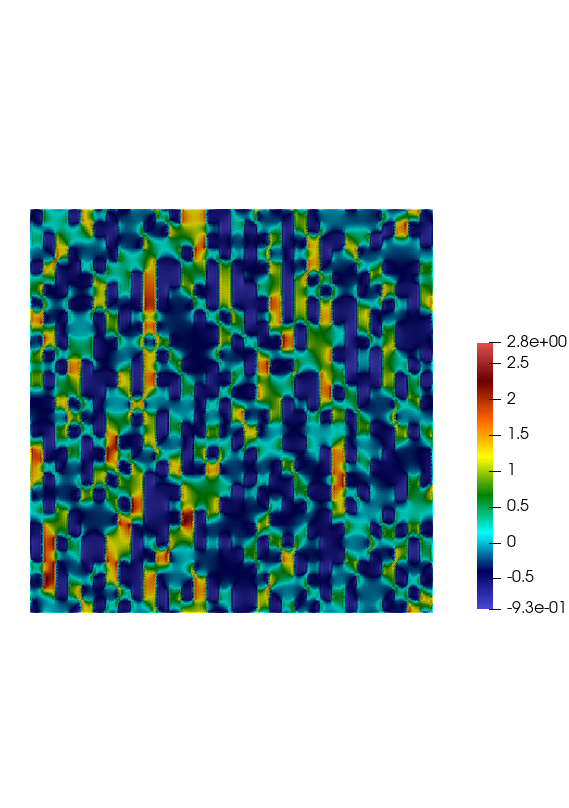}
\quad
\includegraphics[width=0.48\textwidth]{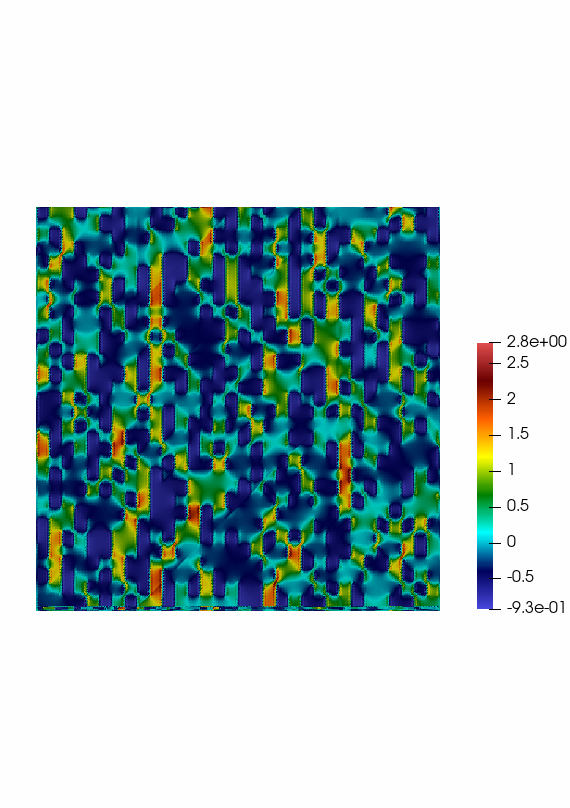}
\caption{Left: approximation of $\partial_{x_1} w_1$ computed using a classical homogenization approach ($\eps = 1/128$, $h=1/256$). Right: approximation $\widetilde{W}^1_\eps$ defined by~\eqref{eq:maison} ($\eps = 1/128$, $h=1/256$, $H = 1/4$). \label{fig:corrector}}
\bigskip
\end{figure}

We now consider a more general case, namely the random checkerboard case. Neither of the two correctors vanish (in contrast to the lamellar case), and both are random functions. For a standard heterogeneous problem of the form $-\dive (k_\eps \nabla \widetilde{u}_\eps) = f$, the two-scale expansion now reads (compare with~\eqref{eq:tse_1correc})
\begin{equation} \label{eq:tse}
\widetilde{u}_\eps(x_1, x_2,\omega) \approx u^\star(x_1, x_2) + \eps \sum_{i=1}^2 w_i(x_1/\eps,x_2/\eps,\omega) \, \partial_{x_i} u^\star(x_1, x_2).
\end{equation}
We follow the approach outlined in the lamellar case and again replace $(\widetilde{u}_\eps, u^\star)$ with the solution $(u_\eps, \overline{u})$ to~\eqref{eq:pb_min_H}. Note that, in view of the boundary conditions $\overline{u}(x) = x_1$ on $\Gamma$, we expect $\overline{u}$ to not depend on $x_2$, and thus the third term in the right-hand side of~\eqref{eq:tse} to vanish. In the vein of~\eqref{eq:corrector_approximated}, we write
\begin{align}
  \partial_{x_1} w_1\left(\frac{x_1}{\eps},\frac{x_2}{\eps},\omega \right) \approx \widetilde{W}^1_\eps(x_1,x_2,\omega) := \frac{\partial_{x_1} u_\eps(x_1,x_2,\omega) - \partial_{x_1} \overline{u}(x_1,x_2)}{\partial_{x_1} \overline{u}(x_1,x_2)},
  \label{eq:maison}
  \\
  \partial_{x_2} w_1\left(\frac{x_1}{\eps},\frac{x_2}{\eps},\omega \right) \approx \widetilde{W}^2_\eps(x_1,x_2,\omega) := \frac{\partial_{x_2} u_\eps(x_1,x_2,\omega) - \partial_{x_2} \overline{u}(x_1,x_2)}{\partial_{x_1} \overline{u}(x_1,x_2)}.
  \nonumber
\end{align}
On Figure~\ref{fig:corrector}, we show, for one particular realization $k_\eps(\cdot,\omega_0)$, the approximation of $\partial_{x_1} w_1(x_1/\eps,x_2/\eps,\omega_0)$ computed by a classical homogenization approach (solving the corrector equation on a large supercell) and the approximation $\widetilde{W}^1_\eps(x_1,x_2,\omega_0)$ defined by~\eqref{eq:maison}. For this case, $k_\eps(\cdot,\omega_0)$ is a realization of the random checkerboard with $\eps = 1/128$ and a contrast $\beta = 9$ (while $k^\star$ has been approximated by the procedure described in Section~\ref{sec:random_nous}, where we have considered $M=16$ realizations of the random checkerboard with $\eps = 1/128$ and $\beta = 9$, on the meshes of size $h=1/256$ and $H=1$; we have found $(\overline{k}^{\rm opt})_{11} = 3.003$ and $(\overline{k}^{\rm opt})_{21} = 0.00736$; for the sake of simplicity, and because it does not change the spirit of our approach, we have fixed $\overline{k}_{12} = k^\star_{12} = 0$ and $\overline{k}_{22} = k^\star_{22} = \sqrt{\beta} = 3$). We qualitatively observe on Figure~\ref{fig:corrector} a good agreement between the two approximations. The relative error, defined by $\dps \frac{\left\| \widetilde{W}_\eps(\cdot,\omega_0) - \partial_1 w_1^h(\cdot/\eps,\omega_0) \right\|_{L^2(D_c)}}{\left\| \partial_1 w_1^h(\cdot/\eps,\omega_0) \right\|_{L^2(D_c)}}$, is equal to 28\%. 
This confirms the interest of~\eqref{eq:maison}, which, at no additional cost, already provides a fairly good approximation of the corrector. We also note that our approach does not rely on any geometrical assumption (such as periodicity, random stationarity, \dots) on the oscillatory coefficient $k_\eps$, in sharp contrast to classical approaches.

The formula~\eqref{eq:maison} hence leads to a fairly accurate approximation of the corrector. We note that a similar accuracy was obtained in~\cite[Sect.~3.3 and Tableau~1]{cras_kun_li}, where we also built an approximation of the corrector using a method without any additional cost. We also point out that, should need be, a better approximation of the corrector could probably be obtained by using a method similar to that used in~\cite[Sect.~3.3]{cocv}: the error is there of the order of a few percents, but the procedure (see Eq.~(3.21) there) is much more expensive than that used here.

\bigskip

\section{Variants and extensions} \label{sec:variants}

This section discusses several variants of the approach described above, concerning the weighting functions in the case when the two models are coupled over a region $D_c$, and concerning an alternative way of coupling the two models. To begin with, and as pointed out above, we note that, in~\eqref{eq:pb_min}, Dirichlet boundary conditions could be imposed on some other part of the exterior boundary of $D$. Likewise, in principle, alternative choices for the boundary conditions on $\Gamma$ (that we have considered to be $\overline{u}(x) = x_1$ on $\Gamma$) could be made. The whole point is to be able to identify analytically the solution to the homogenized problem~\eqref{eq:pb_min_star} in the case $\overline{k} = k^\star$, in order to introduce a relevant $\overline{u}_{\rm ref}$ in~\eqref{eq:def_J}.

\subsection{Weighting functions} \label{sec:variants_weighting}

In Sections~\ref{sec:state_art} and~\ref{sec:nous} above, we have considered an Arlequin coupling between the heterogeneous model with coefficient $k_\eps$ and the tentative effective model with coefficient $\overline{k}$ with the weighting functions $\alpha_1$ and $\alpha_2$ defined by~\eqref{eq:def_alpha}. The energy that we manipulate is thus defined by~\eqref{eq:def_E_pre}.

An alternative choice consists in considering a {\em continuous} weighting function $\alpha_1$ satisfying $\alpha_1 = 1$ in $D$, $\alpha_1 = 0$ in $D_f$ and $0 < \alpha_1(x) < 1$ in $D_c$ (and of course setting $\alpha_2 = 1 - \alpha_1$). 
The qualitative theoretical conclusions drawn above carry over to that case. In particular, we can again optimize upon $\overline{k}$ (in the spirit of~\eqref{eq:optim_J}--\eqref{eq:def_J}) in order to identify $k^\star$.

The algorithmic improvements described in Section~\ref{sec:nous} also carry over to that case. In particular, we can again enrich the Lagrange multiplier space, as explained in Section~\ref{sec:enrich}. Indeed, assuming $\overline{k} = k^\star$ and considering the limit $\eps \to 0$, we observe that the Lagrange multiplier satisfies $\psi = k^\star \, \psi_0$ where $\psi_0$ is independent of $k^\star$ and satisfies (compare to~\eqref{eq:def_psi0})
$$
- \Delta \psi_0 + \psi_0 = \dive (\alpha_1 \, e_1) \ \ \text{in $D_c$},
\qquad
\nabla \psi_0 \cdot n = 0 \ \ \text{on $\partial D_c$}.
$$
As in the case considered in Section~\ref{sec:enrich}, $\psi_0$ can be computed beforehand.

\subsection{Interface coupling} \label{sec:interface}

As mentioned in the Introduction, the heterogeneous model with coefficient $k_\eps$ and the tentative effective model with coefficient $\overline{k}$ can actually be coupled in various ways. The Arlequin coupling technique used in~\cite{cottereau}, which we have described above and which is based on the existence of a subdomain $D_c$ where the two models exist (see Figure~\ref{fig:decompo_left}), is only one possibility among many others. Another possibility is for instance to couple the models through transmission conditions across an interface. In this case, the material is modelled by the oscillatory coefficient $k_\eps$ in the central region $D_f$ that is surrounded by a second region $D$ where the material is modelled by a constant coefficient $\overline{k}$. We then define the energy (compare to~\eqref{eq:def_E_pre}) by
$$
{\cal E}_{\rm interface}(\overline{u},u_\eps) = \frac{1}{2} \int_D \overline{k} \, |\nabla \overline{u}(x)|^2 + \frac{1}{2} \int_{D_f} k_\eps(x) \, |\nabla u_\eps(x)|^2
$$
and consider the minimization problem (compare to~\eqref{eq:pb_min})
\begin{equation} \label{eq:pb_min_interface}
\inf \left\{ \begin{array}{c} {\cal E}_{\rm interface}(\overline{u},u_\eps), \quad \overline{u} \in H^1(D), \quad \overline{u}(x) = x_1 \ \text{on $\Gamma$}, \\ u_\eps \in H^1(D_f), \quad \overline{u} = u_\eps \ \text{on $\partial D_f = \overline{D} \cap \overline{D_f}$} \end{array} \right\},
\end{equation}
where $\Gamma$ is again the same part of the exterior boundary of $D$ as on Figure~\ref{fig:decompo_left}. The two models are now coupled by the condition $\overline{u} = u_\eps$ on the interface $\partial D_f = \overline{D} \cap \overline{D_f}$.

In the limit $\eps \to 0$, Problem~\eqref{eq:pb_min_interface} is well approximated by its homogenized limit
\begin{equation} \label{eq:pb_min_interface_star}
\inf \left\{ \begin{array}{c} {\cal E}_{\rm interface}^\star(\overline{u},u_0), \quad \overline{u} \in H^1(D), \quad \overline{u}(x) = x_1 \ \text{on $\Gamma$}, \\ u_0 \in H^1(D_f), \quad \overline{u} = u_0 \ \text{on $\partial D_f = \overline{D} \cap \overline{D_f}$} \end{array} \right\},
\end{equation}
where the energy ${\cal E}_{\rm interface}^\star$ reads as
$$
{\cal E}_{\rm interface}^\star(\overline{u},u_0) = \frac{1}{2} \int_D \overline{k} \, |\nabla \overline{u}(x)|^2 + \frac{1}{2} \int_{D_f} k^\star \, |\nabla u_0(x)|^2.
$$
We then have the same type of result for this model as for~\eqref{eq:pb_min_star} (see Lemma~\ref{lem:consistence}): if $\overline{k} = k^\star$, the solution to~\eqref{eq:pb_min_interface_star} is $\overline{u}(x) = x_1$ in $D$ and $u_0(x) = x_1$ in $D_f$; conversely, if $(\overline{u},u_0)$ is a solution to~\eqref{eq:pb_min_interface_star} with $\overline{u}(x) = x_1$ in $D$, then $u_0(x) = x_1$ in $D_f$ and $\overline{k} = k^\star$ (see~\cite{ref_olga_theo} for details). 
This thus motivates the idea of again optimizing upon $\overline{k}$ by considering the minimizing problem~\eqref{eq:optim_J} where the functional $J$ is now defined by
\begin{equation} \label{eq:def_J_interface}
  J(\overline{k}) = \int_D | \nabla \overline{u} - e_1 |^2.
\end{equation}
This discussion shows that coupling the heterogeneous model and the constant coefficient model by the Arlequin approach is only one possibility among others, and that, in principle, other choices can be made.

\medskip

We now briefly discuss the discretization of~\eqref{eq:pb_min_interface}. The transmission condition across the interface requires some careful attention. A possibility is to introduce a coarse mesh in $D$ and a fine mesh in $D_f$ which are such that, on the interface $\partial D_f$, the fine mesh is a submesh of the coarse mesh. The transmission condition is then e.g. replaced by the constraint that $\dps \int_{\partial D_f} (\overline{u} - u_\eps) \, \phi = 0$ for any $\phi$ in a suitably chosen finite dimensional space. In the spirit of the enrichment presented in Section~\ref{sec:enrich}, an interesting question is to again choose here a discretization space for the Lagrange multiplier which is well adapted to the problem at hand. 

\medskip

In the same spirit as in Section~\ref{sec:corrector}, it is also possible to design a strategy, in this case of interface coupling, to recover the corrector. In order to gain some intuition, we temporarily consider the one-dimensional case. The solution to~\eqref{eq:pb_min_interface} (for a given $\overline{k}$ and at a given finite value of $\eps$) can then be computed. We observe that, on the interface where the two models co-exist, we have
\begin{equation} \label{eq:ratio}
\frac{u'_\eps(x) - \overline{u}'(x)}{\overline{u}'(x)} = \frac{\overline{k}}{k_\eps(x)} - 1.
\end{equation}
For a fixed $\eps$, the optimal value of $\overline{k}$ (i.e. that which minimizes~\eqref{eq:def_J_interface} where $(\overline{u},u_\eps)$ is the solution to~\eqref{eq:pb_min_interface}) is given by $\dps \overline{k}^{\rm opt}(\eps) = \frac{1}{\langle k_\eps^{-1} \rangle}$, where $\langle k_\eps^{-1} \rangle$ is the average of $k_\eps^{-1}$ on the region $D_f$ where the heterogeneous model is used. We thus obtain
\begin{equation} \label{eq:conv_interface}
  \lim_{\eps \to 0} \overline{k}^{\rm opt}(\eps) = k^\star,
\end{equation}
meaning that the approach indeed computes a converging approximation (when $\eps \to 0$) of the homogenized coefficient. Furthermore, we deduce from~\eqref{eq:ratio} and~\eqref{eq:conv_interface} that, when we consider the coupled problem~\eqref{eq:pb_min_interface} with the optimal value $\overline{k} = \overline{k}^{\rm opt}(\eps)$ and when $\eps \ll 1$, we get
$$
\frac{u'_\eps(x) - \overline{u}'(x)}{\overline{u}'(x)} = \frac{\overline{k}^{\rm opt}(\eps)}{k_\eps(x)} - 1 \approx \frac{k^\star}{k_\eps(x)} - 1 = w'\left(\frac{x}{\eps}\right),
$$
where $w$ is the corrector function of the homogenization theory. We are thus able to recover the corrector by post-processing the solution $(\overline{u},u_\eps)$ of the coupled problem~\eqref{eq:pb_min_interface}. Of course, such a strategy can be extended to problems beyond the one-dimensional case.

\section*{Acknowledgements}

The work of the authors is partially supported by the EOARD under Grant FA9550-17-1-0294. The last two authors are grateful to R\'egis Cottereau for introducing them to the approach of~\cite{cottereau} and for several interesting discussions.

\appendix

\section{Proof of Lemma~\ref{lem:consistence}}
\label{sec:proof_lemma_consistence}
The Euler-Lagrange equations of~\eqref{eq:pb_min_star} read as follows: find $\overline{u} \in V^{\rm Dir BC}$, $u_0 \in H^1(D_c \cup D_f)$ and $\psi \in H^1(D_c)$ such that
\begin{equation} \label{eq:arlequin_homog}
  \begin{cases}
  \forall \overline{v} \in V^0, \ \ & \dps \int_D \overline{k} \, \nabla \overline{u} \cdot \nabla \overline{v} + \int_{D_c} \frac{\overline{k}}{2} \, \nabla \overline{u} \cdot \nabla \overline{v} + C(\overline{v},\psi) = 0,
  \\
  \forall v_0 \in H^1(D_c \cup D_f), \ \ & \dps \int_{D_c} \frac{k^\star}{2} \, \nabla u_0 \cdot \nabla v_0 + \int_{D_f} k^\star \, \nabla u_0 \cdot \nabla v_0 - C(v_0,\psi) = 0,
  \\
  \forall \phi \in H^1(D_c), \ \ & C(\overline{u}-u_0,\phi) = 0,
  \end{cases}
\end{equation}
where $C$ is defined by~\eqref{eq:def_C} and where
\begin{gather*}
V^{\rm Dir BC} = \left\{ v \in H^1(D \cup D_c), \quad \text{$v(x)=x_1$ on $\Gamma$} \right\},
\\
V^0 = \left\{ v \in H^1(D \cup D_c), \quad \text{$v=0$ on $\Gamma$} \right\}.
\end{gather*}
The unique minimizer of~\eqref{eq:pb_min_star} solves~\eqref{eq:arlequin_homog}. Conversely, system~\eqref{eq:arlequin_homog} has a unique solution.

\medskip

We show the first assertion of Lemma~\ref{lem:consistence}. Assume $\overline{k} = k^\star$. Then we check that $(\overline{u},u_0,\psi)$ given by
$$
\overline{u}(x) = x_1 \ \ \text{in $D \cup D_c$}, \qquad u_0(x) = x_1 \ \ \text{in $D_c \cup D_f$},
$$
and
\begin{equation} \label{eq:magic_psi}
\hspace{-3mm} \left\{ \begin{array}{l}
  - \Delta \psi + \psi = 0 \ \text{in $D_c$},
  \\
  \dps \nabla \psi \cdot n = \frac{\overline{k}}{2} \, e_1 \cdot n \ \text{on $\partial D \cap \partial D_c$}, \ \ \nabla \psi \cdot n = -\frac{\overline{k}}{2} \, e_1 \cdot n \ \text{on $\partial D_c \cap \partial D_f$},
\end{array} \right.
\end{equation}
is a solution to~\eqref{eq:arlequin_homog}. Since~\eqref{eq:arlequin_homog} has a unique solution, this implies that $(\overline{u},u_0)$ is the minimizer of~\eqref{eq:pb_min_star}.

\medskip

We now turn to the second assertion of Lemma~\ref{lem:consistence}. Let $(\overline{u},u_0)$ be the unique minimizer of~\eqref{eq:pb_min_star} and assume that $\overline{u}(x) = x_1$ in $D \cup D_c$. Since the minimizer satisfies the Euler-Lagrange equations of the problem, there exists $\psi \in H^1(D_c)$ such that~\eqref{eq:arlequin_homog} holds. The first line of that system reads
$$
\forall \overline{v} \in V^0, \quad \int_D \overline{k} \, e_1 \cdot \nabla \overline{v} + \int_{D_c} \frac{\overline{k}}{2} \, e_1 \cdot \nabla \overline{v} + C(\overline{v},\psi) = 0,
$$
which implies that $\psi$ satisfies~\eqref{eq:magic_psi}.

The third line of~\eqref{eq:arlequin_homog} implies that $u_0(x) = x_1$ in $D_c$, and the second line reads as:
$$
\forall v_0 \in H^1(D_c \cup D_f), \quad \int_{D_c} \frac{k^\star}{2} \, e_1 \cdot \nabla v_0 + \int_{D_f} k^\star \, \nabla u_0 \cdot \nabla v_0 - C(v_0,\psi) = 0.
$$
Using that $\psi$ satisfies~\eqref{eq:magic_psi}, we deduce
\begin{equation} \label{eq:moto}
\begin{cases}
  - \dive (k^\star \nabla u_0) = 0 \quad \text{in $D_f$},
  \\
  \overline{k} \, e_1 \cdot n = k^\star \, e_1 \cdot n \quad \text{on $\partial D \cap \partial D_c$},
  \\
  k^\star \, e_1 \cdot n - 2 k^\star \, \nabla u_0 \cdot n + \overline{k} \, e_1 \cdot n = 0 \quad \text{on $\partial D_c \cap \partial D_f$}.
\end{cases}
\end{equation}
There exists some $x \in \partial D \cap \partial D_c$ such that $e_1 \cdot n(x) \neq 0$. We thus deduce from the second line of~\eqref{eq:moto} that $\overline{k} = k^\star$. The first and third lines of~\eqref{eq:moto} then imply that $u_0(x) = x_1$ in $D_f$. 

\section{Proof of Lemma~\ref{lem:dimanche}} \label{app:equivalence}
We first show that any solution to~\eqref{eq:arlequin_original} is a solution to~\eqref{eq:arlequin_alternative}. To that aim, let $(\overline{u},u_\eps,\psi,\theta)$ be a solution to~\eqref{eq:arlequin_original}. In the second line of~\eqref{eq:arlequin_original}, we choose $v_\eps(x,\omega) = \xi(\omega)$ for some $\xi \in L^2(\Omega)$. Using that $v_\eps$ is independent of $x$, we obtain
$$
\EE[\xi] \int_{D_c} \psi + | D_c | \, \EE[\theta \, \xi] = 0.
$$ 
Since the mean of $\psi$ vanishes, this implies that $\EE[\theta \, \xi] = 0$ for any $\xi \in L^2(\Omega)$ and therefore $\theta = 0$. The first line of~\eqref{eq:arlequin_alternative} is hence satisfied. 

Take now the function $v_\eps(x,\omega) = \xi(\omega) \, \chi(x)$ in the second line of~\eqref{eq:arlequin_original}, for any $\xi \in L^2(\Omega)$ and $\chi \in H^1(D_c \cup D_f)$. Since $\theta = 0$, we have
$$
A_\eps(u_\eps,v_\eps) - C(\EE[v_\eps],\psi) = 0
$$
and hence
$$
\EE \left[ \xi \, A_\eps^\omega(u_\eps,\chi)\right] - \EE\left[ \xi \, C(\chi,\psi) \right] = 0.
$$
This holds for any $\xi \in L^2(\Omega)$: we hence see that the second line of~\eqref{eq:arlequin_alternative} holds almost surely in $\omega$.

Consider next some $\widetilde{\phi} \in H^1(D_c)$ and set $\dps \tau = \frac{1}{|D_c|} \int_{D_c} \widetilde{\phi}$. Then $\phi := \widetilde{\phi} - \tau \in H^1_{\rm mean}(D_c)$ and thus, using the third equation of~\eqref{eq:arlequin_original}, we have $C(\overline{u} - \EE \left[u_\eps \right],\phi) = 0$. We compute
$$
C(\overline{u} - \EE[u_\eps],\widetilde{\phi}) = C(\overline{u} - \EE \left[u_\eps \right],\tau) = \tau \int_{D_c}\left(\overline{u} -\EE \left[u_\eps\right]\right) = 0,
$$
where the last equality stems from the fourth line of~\eqref{eq:arlequin_original} with $\xi(\omega) = \tau$. The third line of~\eqref{eq:arlequin_alternative} is hence satisfied. 

\medskip

We are now left with investigating the uniqueness of the solution to~\eqref{eq:arlequin_alternative}. To that aim, consider two solutions to this problem, denoted $(\overline{u}^1,u^1_\eps,\psi^1)$ and $(\overline{u}^2,u^2_\eps,\psi^2)$. We set $\overline{\bm{u}} = \overline{u}^1 - \overline{u}^2$, $\bm{u}_\eps = u^1_\eps - u^2_\eps$ and $\bm{\psi} = \psi^1 - \psi^2$ and note that $(\overline{\bm{u}},{\bm{u}}_\eps,\bm{\psi})$ satisfies~\eqref{eq:arlequin_alternative} and that $\overline{\bm{u}} \in V^0$, ${\bm{u}}_\eps \in V_f$ and $\bm{\psi} \in H^1(D_c)$. Upon setting $v_\eps = {\bm{u}}_\eps(\cdot,\omega)$ in the second line of~\eqref{eq:arlequin_alternative} and taking the expectation, we deduce
$$ 
A_\eps(\bm{u}_\eps,\bm{u}_\eps) - C(\EE[ \bm{u}_\eps], \bm{\psi}) = 0.
$$
Setting $\overline{v} = \overline{\bm{u}}$ in the first line of~\eqref{eq:arlequin_alternative}, we have
$$ 
\overline{A}(\overline{\bm{u}},\overline{\bm{u}}) + C(\overline{\bm{u}},\bm{\psi}) = 0.
$$
Adding the above two equations and using the third line of~\eqref{eq:arlequin_alternative} yields that $\overline{A}(\overline{\bm{u}},\overline{\bm{u}}) + A_\eps(\bm{u}_\eps, \bm{u}_\eps) = 0$. This implies $\nabla \overline{\bm{u}} = 0$ and $\nabla \bm{u}_\eps = 0$, and hence $\overline{u}^1 = \overline{u}^2$ and $\bm{u}_\eps(\cdot,\omega) = c_1(\omega)$ almost surely. The third line of~\eqref{eq:arlequin_alternative} now yields that $C(\EE[c_1(\omega)],\phi) = 0$ for any $\phi \in H^1(D_c)$, which implies $\EE[c_1(\omega)] = 0$.

Consider eventually the second line of~\eqref{eq:arlequin_alternative} with $v_\eps = \bm{\psi}$. Since $\nabla \bm{u}_\eps = 0$, we have $C(\bm{\psi},\bm{\psi}) = 0$ and hence $\psi^1 = \psi^2$.

We have hence shown the uniqueness of a solution to~\eqref{eq:arlequin_alternative} up to the addition, for $u_\eps$, of a random constant of vanishing expectation.

\section{Initial guess formula in the case of matrix-valued coefficients} \label{sec:2Dformula}
We present here a formula analogous to~\eqref{eq:IG_k1} for the case when the constant coefficient upon which we optimize is matrix-valued: $\overline{k} = \begin{bmatrix} \overline{k}_{11} & \overline{k}_{12} \\ \overline{k}_{21} & \overline{k}_{22} \end{bmatrix}$. Similarly to the scalar-case, the Arlequin equations amount to solving~\eqref{eq:arlequin_pena_matrix} where $\overline{\cal A}$ depends linearly on $\overline{k}$ while $\overline{F}$, $\overline{\cal P}$, ${\cal A}_\eps$, $\overline{\cal C}$ and ${\cal C}_\eps$ are independent of $\overline{k}$. We thus write~\eqref{eq:arlequin_pena_matrix} in the form (compare to~\eqref{eq:sys_lin})
\begin{equation} \label{eq:sys_lin_2D}
\Big( \overline{k}_{11} \, Z_1 + \overline{k}_{12} \, Z_2 + \overline{k}_{21} \, Z_3 + \overline{k}_{22} \, Z_4 + Z_5 \Big) \, U(\overline{k}) = F_0,
\end{equation}
where the vector $F_0$ and the matrices $Z_1$, $Z_2$, $Z_3$, $Z_4$ and $Z_5$ do not depend on $\overline{k}$. Let $U(\overline{k}^0)$ be the solution to~\eqref{eq:sys_lin_2D} for some matrix $\overline{k}^0$. Proceeding as in~\eqref{eq:moto2}, we have
$$
M(\overline{k}) \, U(\overline{k}) = U(\overline{k}^0),
$$
where the matrix $M(\overline{k})$ depends on $\overline{k}$ in an affine manner. Setting $A_{\overline{k}} = \overline{k}_{11} Z_1 + \overline{k}_{12} Z_2 + \overline{k}_{21} Z_3 + \overline{k}_{22} Z_4 + Z_5$ for any matrix $\overline{k}$, we have $\dps M(\overline{k}) = \left( A_{\overline{k}^0} \right)^{-1} A_{\overline{k}}$ and hence (compare with~\eqref{eq:def_M})
\begin{multline*}
  M(\overline{k}) =
  (\overline{k}_{11} - \overline{k}^0_{11}) \left( A_{\overline{k}^0} \right)^{-1} Z_1 +
  (\overline{k}_{12} - \overline{k}^0_{12}) \left( A_{\overline{k}^0} \right)^{-1} Z_2
  \\
  + (\overline{k}_{21} - \overline{k}^0_{21}) \left( A_{\overline{k}^0} \right)^{-1} Z_3 +
  (\overline{k}_{22} - \overline{k}^0_{22}) \left( A_{\overline{k}^0} \right)^{-1} Z_4 + \Id,
\end{multline*}
where $\Id$ is the identity matrix.

\medskip

As in Section~\ref{sec:IG}, we would like to find $\overline{k}$ such that $U(\overline{k}) = U_{\rm target}$ for some given $U_{\rm target}$. We again recast the problem as finding $\overline{k}$ that minimizes the function $f(\overline{k}) = \| U(\overline{k}_0) - M(\overline{k}) \, U_{\rm target} \|^2$. We thus ask for the partial derivatives of $f$ to vanish and we obtain that $\overline{k}$ should satisfy the following linear system:
\begin{equation}
 \label{linsyst_2D}
 \begin{cases}
   C_{11} \, \overline{k}_{11} + C_{12} \, \overline{k}_{12} + C_{13} \, \overline{k}_{21} + C_{14} \, \overline{k}_{22} = B_1,
   \\
   C_{21} \, \overline{k}_{11} + C_{22} \, \overline{k}_{12} + C_{23} \, \overline{k}_{21} + C_{24} \, \overline{k}_{22} = B_2,
   \\
   C_{31} \, \overline{k}_{11} + C_{32} \, \overline{k}_{12} + C_{33} \, \overline{k}_{21} + C_{34} \, \overline{k}_{22} = B_3,
   \\
   C_{41} \, \overline{k}_{11} + C_{42} \, \overline{k}_{12} + C_{43} \, \overline{k}_{21} + C_{44} \, \overline{k}_{22} = B_4,
 \end{cases}
\end{equation}
where $\dps C_{ij} = - \left( A_{\overline{k}^0} \right)^{-1} Z_i \, U_{\rm target} \cdot \left( A_{\overline{k}^0} \right)^{-1} Z_j \, U_{\rm target}$ for all $1 \leq i,j \leq 4$ and where $\dps B_i = \left( A_{\overline{k}^0} \right)^{-1} Z_i \, U_{\rm target} \cdot \left[ \left( A_{\overline{k}^0} \right)^{-1} Z_5 \, U_{\rm target} - U(\overline{k}^0) \right]$ for any $1 \leq i \leq 4$.

As in Section~\ref{sec:IG}, we eventually replace the vectors $U_{\rm target}$ and $U(\overline{k}^0)$ by $\pi \left(U_{\rm target} \right)$ and $\pi \left(U(\overline{k}^0) \right)$ in~\eqref{linsyst_2D}. The resulting solution $(\overline{k}_{11},\overline{k}_{12},\overline{k}_{21},\overline{k}_{22})$ of the linear system is then considered as an adequate initial guess for the optimization problem~\eqref{eq:optim_J}.


\bibliographystyle{plain}
\bibliography{biblio_standard_layout}

\end{document}